\documentclass[preprints,article,accept,moreauthors]{mdpi}

\newcommand{\dif}[0]{{\rm d}}

\DeclareMathOperator{\Law}{Law}

\firstpage{1} 
\pubvolume{xx}
\issuenum{1}
\articlenumber{5}
\pubyear{2019}
\copyrightyear{2019}
\history{Received: date; Accepted: date; Published: date}

\Title{State and Parameter Estimation from Observed Signal Increments}

\Author{Nikolas N\"usken $^{1}$, Sebastian  Reich $^{1}$* and Paul J.\ Rozdeba  $^{1}$}

\AuthorNames{Nikolas N\"usken, Sebastian Reich and Paul J.\ Rozdeba}

\address{%
$^{1}$ \quad Institute of Mathematics, University of Potsdam, Karl-Liebknecht-Str. 24/25, D-14476 Potsdam, Germany}

\corres{Correspondence: sereich@uni-potsdam.de; Tel.: +49-331-977-1859}

\abstract{The success of the ensemble Kalman filter has triggered a strong interest
in expanding its scope beyond classical state estimation problems. In this paper, we focus
on continuous-time data assimilation where the model and measurement errors are correlated
and both states and parameters need to be identified. Such scenarios arise from noisy and partial 
observations of Lagrangian particles which move under a stochastic velocity field involving unknown
parameters. We take an appropriate class of McKean--Vlasov equations as the starting 
point to derive ensemble Kalman--Bucy filter algorithms for combined state and parameter estimation. 
We demonstrate their performance through a series of increasingly complex multi-scale model systems. }


\begin{document}


\section{Introduction}
The research presented in this paper has been motivated by the state and parameter estimation problem for particles moving under a stochastic
velocity field, with the measurements given by partial and noisy observations of their position increments. If
the deterministic contributions to the velocity field are stationary, and the position increments of the moving
particle are exactly observed, then one is led to a standard parameter estimation problem for stochastic differential equations 
(SDEs) \citep{sr:Kut04,sr:P14}. In \cite{sr:AHSV07}, this setting was extended to the case where the deterministic 
contributions to the velocity field themselves undergo a stochastic time evolution. Furthermore, while continuous-time 
observations of position increments are at the focus of the present study, the assimilation of discrete-time observations 
of particle positions has been investigated in \cite{sr:SKJI06,sr:AJS08} under a so-called Lagrangian data assimilation 
setting for atmospheric fluid dynamics.

The assumption of exactly and fully observed position increments is not always realistic and the case of partial and
noisy observations is at the center of the present study. Having access to partial and noisy observations of position increments 
leads to correlations between the measurement and model errors. The theoretical impact of such correlations on state and parameter
estimation problems has been discussed, for example, in \cite{sr:simon} in the context of linear systems, and in \cite{sr:BC09} for nonlinear systems. 
One finds in particular that the appropriately adjusted data likelihood involves the gradient of log-densities,
which is nontrivial from a computational perspective, and which prevents a straightforward application
of standard Markov chain Monte Carlo (MCMC) or sequential Monte Carlo (SMC) methods \cite{sr:Liu}.

In this paper, we instead follow an alternative Monte Carlo approach based on appropriately adjusted 
McKean--Vlasov filtering equations, an approach pioneered in \cite{sr:crisan10} in the context of the standard state
estimation problem for diffusion processes. We recall that the notion of McKean--Vlasov equations, first studied in \cite{sr:McK66}, 
characterises a class of SDEs for which their right-hand side depends on the law of the process itself.
We rely on a particular formulation of such McKean--Vlasov filtering equations, the so-called feedback particle filters \citep{sr:meyn13}, 
utilising stochastic innovation processes \citep{sr:reich19}. Our proposed Monte Carlo formulation avoids the need for estimating log-densities, 
and can be implemented in a numerically robust manner relying on a generalised ensemble Kalman--Bucy filter approximation  
applied to an extended state space formulation \citep{sr:majda}. The ensemble Kalman--Bucy filter \citep{sr:br11,sr:TdWMR17}
has been introduced previously as an extension of the popular ensemble Kalman filter \cite{sr:majda,sr:stuart15,sr:reichcotter15} 
to continuous-time data assimilation under the assumption of uncorrelated measurement and model errors. 

We apply the proposed algorithms to a series of state and parameter estimation problems of increasing complexity.
First, we study the state and parameter estimation problem for an Ornstein--Uhlenbeck process \cite{sr:P14}. Two
further experiments investigate the behaviour of the filters for reduced model equations, with the data being 
collected from underlying multi-scale models. There we distinguish between the averaging and homogenisation
scenarios \citep{sr:ps}. Finally, we also look at nonparametric drift estimation \cite{sr:AHSV07}, and parameter estimation for
the stochastic heat equation \cite{sr:AR19}.

%
\section{Mathematical problem formulation}
%

We consider the time evolution of a random state variable $X_t \in \mathbb{R}^{N_x}$ in $N_x$-dimensional state space, $N_x \ge 1$,
as prescribed by an SDE of the form
\begin{equation}
\label{eq:dynamics}
{\rm d}X_t = f(X_t,a)\,{\rm d}t + G \,{\rm d} W_t,
\end{equation}
for time $t \ge 0$, with the drift function $f:\mathbb{R}^{N_x}\times \mathbb{R}^{N_a} \to \mathbb{R}^{N_x}$ depending on $N_a\ge 0$ 
unknown parameters $a = (a^1,\ldots,a^{N_a})^{\rm T} \in \mathbb{R}^{N_a}$. Model errors are represented through 
standard $N_w$-dimensional Brownian motion $W_t$,
$N_w \ge 1$, and a matrix $G \in \mathbb{R}^{N_x\times N_w}$. We also introduce the associated model error covariance matrix 
$Q = G G^{\rm T}$. We will generally assume that the initial condition $X_0$ is fixed, that is, $X_0 = x_0$ a.s.~for given $x_0 \in \mathbb{R}^{N_x}$. 
In terms of a more specific example, one can think of $X_t$ denoting the position of a particle at time $t\ge 0$ moving in $N_x=3$ dimensional space 
under the influence of a stochastic velocity field, with deterministic contributions given by $f$ and stochastic perturbations by $GW_t$.  In the case $G = 0$, the SDE (\ref{eq:dynamics}) reduces to 
an ordinary differential equation with given initial condition $x_0$.

We assume throughout this paper that (\ref{eq:dynamics}) possesses unique, strong solutions for all parameter values $a$. 
See, for example, \cite{sr:P14} for necessary conditions on the drift function $f$. The distribution of $X_t$ is denoted by $\pi_t$, which
we also abbreviate by $\pi_t = {\rm Law}(X_t)$.  We use the same notation for measures and their Lebesgue densities, provided they exist.

\begin{Example} \label{ex:example0}
A wide class of drift functions can be written in the form
\begin{equation} \label{eq:linear_parameter}
f(x,a) = f_0(x) + B(x) a = f_0(x) + \sum_{i=1}^{N_a} b_i(x) a^i,
\end{equation}
where $f_0:\mathbb{R}^{N_x} \to \mathbb{R}^{N_x}$ is a known drift function, 
the $b_i:\mathbb{R}^{N_x} \to \mathbb{R}^{N_x}$, $i=1,\ldots,N_a$, denote appropriate basis functions,
and the vector $a = (a^1,\ldots,a^{N_a})^{\rm T} \in \mathbb{R}^{N_a}$ contains the unknown parameters of the model. The family 
$\{b_i(x)\}$ of basis functions, which we collect in a matrix-valued function $B(x) = (b_1(x),b_2(x),\ldots,b_{N_a}(x))\in 
\mathbb{R}^{N_x\times N_a}$, could arise from a finite-dimensional truncation of some appropriate Hilbert space $\mathcal{H}$. See, for example, 
\cite{sr:PPRS12} for computational approaches to nonparametric drift estimation using a Galerkin approximation in $\mathcal{H}$, 
where the $b_i(x)$ become finite element basis functions. Furthermore, the expansion coefficients $\{a^i\}$ could be made
time-dependent by letting them evolve according to some system of differential equations arising, for example, from the discretisation
of an underlying partial differential equation with solutions in $\mathcal{H}$. 
See \cite{sr:AHSV07} for specific examples of such a setting. While the present paper focuses on
stationary drift functions, that is, the parameters $\{a^i\}$ are time-independent, the results from Sections 
\ref{sec:PE} and \ref{sec:SPE}, respectively, can easily be extended to the non-stationary case where the parameters themselves 
satisfy given evolution equations. 
\end{Example}

Data and an observation model are required in order to perform state and parameter estimation for SDEs of the form (\ref{eq:dynamics}).
In this paper, we assume that we observe partial and noisy increments ${\rm d}Y_t$ of the signal $X_t$, that is,
\begin{equation}
\label{eq:obs}
{\rm d}Y_t = H \,{\rm d} X_t + R^{1/2}{\rm d}V_t = Hf(X_t,a)\,{\rm d}t + HG\, {\rm d}W_t + R^{1/2}{\rm d}V_t, \quad Y_0 = X_0 = x_0,
\end{equation}
for $t$ in the observation interval  $[0,T]$, $T>0$, 
where $H \in \mathbb{R}^{N_y \times N_x}$ is a given linear operator, $V_t$ denotes standard $N_y$-dimensional
Brownian motion with $N_y \ge 1$ and $R \in \mathbb{R}^{N_y \times N_y}$ is a covariance matrix. We introduce the observation map
\begin{equation} \label{eq:forward}
h(x,a) = Hf(x,a)
\end{equation}
for later use. Unless $HG = 0$, we find that 
the model error $E_t^{\rm m} := GW_t$ in (\ref{eq:dynamics}) and the total observation error 
\begin{equation} \label{eq:model_error}
E_t^{\rm o} := HG W_t + R^{1/2} V_t
\end{equation}
in (\ref{eq:obs}) are correlated. The impact of correlations between the model and measurement errors on the state estimation problem 
have been discussed by \cite{sr:simon,sr:BC09}. Furthermore, such correlations require adjustments to sequential estimation
methods \citep{sr:sarkka,sr:stuart15,sr:reichcotter15} which are the main focus of this paper. We assume throughout this paper
that the covariance matrix 
\begin{equation} \label{eq:measurement_C}
C = HGG^{\rm T}H^{\rm T} + R = HQH^{\rm T} + R
\end{equation} 
of the observation error (\ref{eq:model_error}) is invertible. 

The special case $R=0$ and $H=I$ leads to a pure parameter estimation problem, which has been extensively studied in the literature in the settings 
of maximum likelihood and Bayesian estimators \citep{sr:Kut04,sr:P14}. We will provide a reformulation of the Bayesian approach in the form 
of McKean--Vlasov equations in the parameters, based on the results in \cite{sr:crisan10,sr:meyn13}  in Section \ref{sec:PE}.

If $R \not=0$, then (\ref{eq:dynamics}) and (\ref{eq:obs}) lead to a combined state and parameter
estimation problem with correlated noise terms. We will first discuss the impact of this correlation on the pure state estimation problem in Section 
\ref{sec:state_estimation} assuming that the parameters of the problem are known. Again, we will derive appropriate McKean--Vlasov equations 
in the state variables. Our key contribution is a formulation that avoids the need for log-density estimates, and can be
put into an appropriately generalised ensemble Kalman--Bucy filter approximation framework  \citep{sr:br11,sr:TdWMR17}.
We also formally demonstrate that the McKean--Vlasov filter equation reduces to ${\rm d}X_t = {\rm d}Y_t$ in the limit $R\to 0$ and $H=I$, a property
which is less straightforward to demonstrate for filter formulations involving log-densities.

These McKean--Vlasov equations can be generalised to the combined state and parameter estimation problem 
via an augmentation of state space \citep{sr:majda} in Section \ref{sec:SPE}. Given the results from Section \ref{sec:state_estimation}, 
such an extension is rather straightforward. 

The numerical experiments in Section \ref{sec:Examples} rely exclusively on the generalised ensemble Kalman--Bucy filter approximation
to the McKean--Vlasov equations, which are easy to implement and yield robust and accurate numerical results.


%
%
\section{Parameter estimation from noiseless data} \label{sec:PE}
%
%

In this section, we treat the simpler Bayesian parameter estimation problem which arises from setting
$R=0$ and $H=I$ in (\ref{eq:obs}), that is, $N_y = N_x$. This leads to ${\rm d}X_t = {\rm d}Y_t$ and, furthermore, $X_t = Y_t$ for all $t \in 
[0,T]$, provided $X_0 = Y_0 = x_0$ which we assume throughout this paper.  The requirement that $C = Q$ is invertible requires
that $G$ has rank $N_x$, that is, $N_w \ge N_x$ in (\ref{eq:dynamics}). The data likelihood
\begin{equation}
l_t(a) = \exp \left( \int_0^t f(Y_s,a)^{\rm T} Q^{-1} \mathrm{d}Y_s - \frac{1}{2} \int_0^t   f(Y_s,a)^{\rm T} Q^{-1} f(Y_s,a)\,{\rm d}s \right)
\end{equation}
thus follows from the observation model with additive Brownian noise in (\ref{eq:obs}). 
Given a prior distribution $\Pi_0(a)$ for the parameters, the resulting posterior distribution at any time $t \in (0,T]$ is
\begin{equation} \label{eq:time-dependent-PDF}
\Pi_t(a) = \frac{l_t(a)\Pi_0(a)}{\Pi_0[l_t]}
\end{equation}
according to Bayes' theorem \cite{sr:BC09}. Here, we have introduced the shorthand
\begin{equation}
\Pi_0[l_t] = \int_{\mathbb{R}^{N_a}} l_t(a)\Pi_0(a)\,{\rm d}a
\end{equation}
for the expectation of $l_t$ with respect to $\Pi_0$. It is well-known that 
the posterior distributions $\Pi_t$ satisfy the stochastic partial differential equation
\begin{equation} \label{eq:KSEa}
{\rm d} \Pi_t[\phi] = \left( \Pi_t[ \phi\,h_t] - \Pi_t[\phi]\Pi_t[h_t]\right)^{\rm T} Q^{-1} ({\rm d}Y_t - \Pi_t[h_t]{\rm d}t)
\end{equation}
with time-dependent observation map 
\begin{equation}
h_t (a) = f(Y_t,a), 
\end{equation}
where $\phi:\mathbb{R}^{N_a}\to \mathbb{R}$ is a compactly supported smooth test function, and $\Pi_t[\phi]$ again 
denoting the expectation of $\phi$ with respect to $\Pi_t$. See \cite{sr:BC09} for a detailed discussion. Equation (\ref{eq:KSEa}) constitutes
a special instance of the well-known Kushner--Stratonovitch equation from time-continuous filtering \cite{sr:BC09}. 

%
%
\subsection{Feedback particle filter}
%
%

We now state a McKean--Vlasov reformulation of the Kushner--Stratonovitch equation (\ref{eq:KSEa}) as a special instance of the 
feedback particle filter of \cite{sr:meyn13,sr:reich19}. The key idea is to formulate a stochastic differential equation in the parameters in which they are treated as time-dependent random variables. We introduce the notation $\widetilde{A}_t$ for these, and require that the law of $\widetilde{A}_t$ coincide with (\ref{eq:time-dependent-PDF}) for $t\in [0,T]$, that is, with the solution to (\ref{eq:KSEa}).

\begin{Lemma}[Feedback particle filter]
	\label{lem:fpf}
	 Consider the McKean--Vlasov equations
\begin{equation} \label{eq:fpf1I}
\mathrm{d} \widetilde{A}_t = K_t (\widetilde{A}_t) \, \mathrm{d}I_t + \Omega_t(\widetilde{A}_t) \, \mathrm{d}t,
\end{equation}
where the matrix-valued Kalman gain $K_t \in \mathbb{R}^{N_a\times N_y}$ satisfies
\begin{equation}
\label{eq:gain_equation}
\nabla \cdot \left( \widetilde{\Pi}_t \left( K_t Q \right)  \right) = - \widetilde{\Pi}_t \left(  h_t  - \widetilde{\Pi}_t[h_t] \right)^{\rm T}, 
\quad \widetilde{\Pi}_t = \mathrm{Law}(\widetilde{A}_t),
\end{equation} 
the innovation process $I_t$ can be chosen to be given by either
\begin{equation} \label{eq:innovation1}
\mathrm{d}I_t = \mathrm{d}Y_t - \frac{1}{2} \left( h_t(\widetilde{A}_t) + \widetilde{\Pi}_t[h_t]  \right) \mathrm{d}t,
\end{equation}
or
\begin{equation} \label{eq:innovation2}
{\rm d}I_t =  \mathrm{d}Y_t - \left\{ h_t(\widetilde{A}_t)\, \mathrm{d}t  + G\, {\rm d}W_t  \right\} ,
\end{equation}
and 
\begin{equation}
\Omega^i_t = \frac{1}{2} \sum_{j=1}^{N_a} \sum_{k,l = 1}^{N_y} Q^{kl} K_t^{jl} \left( \partial_{j} K_t^{ik} \right), \quad i = 1,\ldots, N_a.
\end{equation} 	
Then, the distribution $\widetilde{\Pi}_t = \Law(\widetilde{A}_t)$ coincides with the solution to \eqref{eq:KSEa}, 
provided that the initial distributions agree. In other words, $\widetilde{\Pi}_t = \Pi_t$ for all $t \in [0,T]$.
\end{Lemma}

Throughout this paper, we write (\ref{eq:fpf1I}) in the more compact Stratonovitch form
\begin{equation} \label{eq:fpf1S}
\mathrm{d} \widetilde{A}_t = K_t (\widetilde{A}_t)  \circ \mathrm{d}I_t ,
\end{equation}
where the Stratonovitch interpretation is to be applied only to $\widetilde{A}_t$ in $K_t(\widetilde{A}_t)$, while the explicit time-dependence of
$K_t$ remains in its It\^o interpretation. It should be noted that the matrix-valued function $K_t$ is not uniquely defined by the 
PDE (\ref{eq:gain_equation}).  Indeed, provided $K_t$ solves (\ref{eq:gain_equation}), $K_t + \beta_t$ is also a solution whenever $\nabla \cdot \left(\widetilde{\Pi}_t \beta_t \right) = 0$.  As discussed in \cite{sr:TdWMR17}, the minimiser over all suitable $K_t$ with respect to a kinetic energy-type functional is of the form
\begin{equation}
\label{eq:K Psi}
K_t = \nabla \Psi_t Q^{-1},
\end{equation}
for a vector of potential functions $\Psi_t = (\psi_t^1,\ldots, \psi_t^{N_x})$, $\psi^k_t:\mathbb{R}^{N_a}\to\mathbb{R}$. 
Inserting \eqref{eq:K Psi} into \eqref{eq:gain_equation} leads to $N_x$ elliptic partial differential equations (often referred to as Poisson equations),
\begin{equation}
\label{eq:Poisson}
\nabla \cdot \left( \widetilde{\Pi}_t \nabla \Psi_t \right) = - \widetilde{\Pi}_t\left( h_t - \widetilde{\Pi}_t[h_t]\right)^{\rm T}, \quad \widetilde{\Pi}_t[\Psi_t] = 0,
\end{equation}
understood componentwise,
where the centering condition $\widetilde{\Pi}_t[\Psi_t] = 0$ makes the solution unique under mild assumptions on $\widetilde{\Pi}_t$, see \cite{laugesen2015poisson}.
Finally, (\ref{eq:innovation2}) yields a particularly appealing formulation, since it is based on a direct comparison of ${\rm d}Y_t$ with a 
random realisation of the right hand side of the SDE (\ref{eq:dynamics}), given a parameter value $a = \widetilde{A}_t(\omega)$ and
a realisation of the noise term ${\rm d} W_t(\omega)$. This fact will be explored further  in Section \ref{sec:state_estimation}.

\begin{Remark}
	For clarity, let us repeat equations \eqref{eq:gain_equation} and \eqref{eq:K Psi} in their index forms:
\begin{equation}
\label{eq:gain_equation_index}
\sum_{i=1}^{N_a}\sum_{j=1}^{N_y}\partial_{i} \left( \widetilde{\Pi}_t \left( K_t^{ij}Q^{jk} \right)  \right) = - \widetilde{\Pi}_t \left(  h_t^k  - \widetilde{\Pi}_t[h_t^k] \right), \quad 
k = 1, \ldots, N_y, 
\end{equation} 


\begin{equation}
\sum_{j=1}^{N_y} K_t^{ij}(a)Q^{jk} = \partial_{i} \psi^k_t(a), \quad i=1,\ldots, N_a,\quad k=1,\ldots,N_y.
\end{equation}
\end{Remark}

%
%
\subsection{Ensemble Kalman--Bucy filter}
%
%

Let us now assume that the initial distribution $\Pi_0$ is Gaussian, and that $f$ is linear in the unknown parameters such as in (\ref{eq:linear_parameter}). 
Then, the distributions $\widetilde{\Pi}_t$ remain Gaussian for all times with mean $\overline{a}_t$ and covariance matrix $P^{aa}_t$. The elliptic PDE
(\ref{eq:gain_equation}) is solved by the parameter-independent Kalman gain matrix
\begin{equation} \label{eq:Kalman_gain1}
K_t = P_t^{aa} B(Y_t)^{\rm T} Q^{-1}
\end{equation}
and one obtains the McKean--Vlasov formulation 
\begin{equation} \label{eq:KBMFE1}
{\rm d}\widetilde{A}_t = P_t^{aa} B(Y_t)^{\rm T} Q^{-1} {\rm d} I_t
\end{equation} 
of the Kalman--Bucy filter, with the innovation process $I_t$ defined by either
\begin{equation}
 {\rm d}I_t = {\rm d}Y_t - \left( f_0(Y_t) + \frac{1}{2} B(Y_t)(\widetilde{A}_t+\overline{a}_t)\right) {\rm d}t 
\end{equation}
or
\begin{equation}
{\rm d} I_t = {\rm d}Y_t - \left\{\left( f_0(Y_t) + B(Y_t)\widetilde{A}_t\right){\rm d}t + G\,{\rm d} W_t \right\}.
\end{equation}
Note that the Stratonovitch formulation (\ref{eq:fpf1S}) reduces to the standard It\^o interpretation, since $K_t$ no longer depends explicitly on 
$\widetilde{A}_t$.

The McKean--Vlasov equations (\ref{eq:KBMFE1}) can be extended to nonlinear, non-Gaussian parameter estimation problems by 
generalising the parameter-independent Kalman gain matrix (\ref{eq:Kalman_gain1}) to
\begin{equation} \label{eq:Kalman_gain2}
K_t = P_t^{ah} Q^{-1}, \quad P_t^{ah} = \widetilde{\Pi}_t \left[(a-\overline{a}_t)(h_t(a)-\widetilde{\Pi}_t[h_t])^{\rm T}\right] 
=  \widetilde{\Pi}_t \left[ a\,(h_t(a)-\widetilde{\Pi}_t[h_t])^{\rm T}\right] 
\end{equation}
Clearly, the gain (\ref{eq:Kalman_gain2}) provides only an approximation to the solution of  (\ref{eq:gain_equation}). However, such
approximations have become popular in nonlinear state estimation in the form of the ensemble Kalman filter \citep{sr:stuart15,sr:reichcotter15},
and we will test its suitability for parameter estimation in Section \ref{sec:Examples}.

Numerical implementations of the proposed McKean--Vlasov approaches rely on Monte--Carlo approximations. More specifically, given $M$ samples
$\widetilde{A}_0^i$, $i = 1, \ldots, M$, from the initial distribution $\Pi_0$, we introduce the interacting particle system 
\begin{equation}
{\rm d} \widetilde{A}_t^i = K_t^M(\widetilde{A}_t^i)\circ {\rm d}I_t^i,
\end{equation}
where the innovation processes $I_t^i$ are defined by either
\begin{equation}
{\rm d} I_t^i = {\rm d}Y_t - \frac{1}{2} \left( h_t(\widetilde{A}_t^i) + \overline{h}_t^M\right){\rm d}t, \qquad 
\overline{h}_t^M = \frac{1}{M} \sum_{i=1}^M h_t(\widetilde{A}_t^i),
\end{equation}
or, alternatively,
\begin{equation}
{\rm d} I_t^i = {\rm d}Y_t - \left( h_t(\widetilde{A}_t^i)\,{\rm d}t + G \,{\rm d} W_t^i \right),
\end{equation}
and $W_t^i$, $i=1,\ldots,M$, denote independent $N_w$-dimensional Brownian motions.
For $K^M_t$, we will use the parameter-independent empirical Kalman gain approximation
\begin{equation}
K_t^M =  \widehat{P}_t^{ah}Q^{-1} , \qquad \widehat{P}_t^{ah} = \frac{1}{M-1} \sum_{i=1}^M \widetilde{A}_t^i(h_t(\widetilde{A}_t^i)-\overline{h}_t^M)^{\rm T},
\end{equation}
in our numerical experiments, which leads to the so-called ensemble Kalman--Bucy filter
\citep{sr:br11,sr:TdWMR17}. Note that $\widehat{P}_t^{ah}$ provides an unbiased estimator of $P_t^{ah}$.

Finally, a robust and efficient time-stepping procedure for approximating $\widetilde{A}_{t_n}$, $t_n = n\Delta t$, 
is provided in \citep{sr:akir11,sr:dWRS18,sr:BSW18}. Denoting the approximations at time $t_n$ by $\widetilde{A}_n^i$, $i=1,\ldots,M$, we obtain
\begin{equation} \label{eq:EnKBF1a}
\widetilde{A}_{n+1}^i = \widetilde{A}_n^i + \Delta t \widehat{P}_n^{ah}\left( Q + \Delta t \widehat{P}_n^{hh}\right)^{-1} \Delta I_n^i
\end{equation}
with step size $\Delta t>0$, empirical covariance matrices
\begin{equation}
\widehat{P}_n^{ah} = \frac{1}{M-1} \sum_{i=1}^M \widetilde{A}_n^i(h_n(\widetilde{A}_n^i)-\overline{h}_n^M)^{\rm T}, \qquad
\widehat{P}_n^{hh} = \frac{1}{M-1} \sum_{i=1}^M h_n(\widetilde{A}_n^i)(h_n(\widetilde{A}_n^i)-\overline{h}_n^M)^{\rm T},
\end{equation}
and innovation increments $\Delta I_n^i$ given by either
\begin{equation}
\Delta I_n^i = \Delta Y_n - \frac{1}{2} \left( h_n(\widetilde{A}_n^i) + \overline{h}_n^M\right)\Delta t , \qquad
\overline{h}_n^M = \frac{1}{M} \sum_{i=1}^M h_n(\widetilde{A}_n^i),
\end{equation}
or
\begin{equation}
\Delta  I_n^i = \Delta Y_n - \left( h_n(\widetilde{A}_n^i)\,\Delta t + \Delta t^{1/2} G \Xi_n^i \right), \qquad \Xi_n^i \sim {\rm N}(0,I).
\end{equation}
Here we have used the abbreviations $h_n(a) = f(Y_n,a)$, $Y_n = Y_{t_n}$, and $\Delta Y_n = Y_{t_{n+1}}-Y_{t_n}$.

While the feedback particle formulation (\ref{eq:fpf1S}) and its ensemble Kalman--Bucy filter approximation (\ref{eq:EnKBF1a}) are special
cases of already available formulations, they provide the starting point for our novel McKean--Vlasov equations and their numerical approximation of the 
combined state and parameter estimation problem with correlated measurement and model errors, which we develop in the following two sections.

%
%
\section{State estimation for noisy data} \label{sec:state_estimation}
%
%

We return to the observation model (\ref{eq:obs}) with $R \not=0$ and general $H$. The pure state estimation problem is considered first, 
that is, $f(x,a) = f(x)$ in (\ref{eq:dynamics}).

Using $E_t^{\rm o}$, given by (\ref{eq:model_error}), 
and $E_t^{\rm c}$ defined by
\begin{equation}
E_t^{\rm c} = G(I-G^{\rm T}H^{\rm T}C^{-1} H G) W_t - QH^{\rm T}C^{-1} R^{1/2} V_t
\end{equation}
with the total measurement error covariance matrix $C$ given by (\ref{eq:measurement_C}),
we find that
\begin{equation}
G W_t = E_t^{\rm c} + QH^{\rm T} C^{-1} E_t^{\rm o},
\end{equation}
and the covariations \cite{sr:P14} satisfy
\begin{equation}
\langle E^{\rm o},E^{\rm c}\rangle_t = 0, \quad \langle E^{\rm o},E^{\rm o}\rangle_t = Ct, \quad
\langle E^{\rm c},E^{\rm c}\rangle_t = G(I-G^{\rm T}H^{\rm T}C^{-1} H G)G^{\rm T}t.
\end{equation}
Hence (\ref{eq:dynamics}) and (\ref{eq:obs}) can be rewritten as follows:
\begin{subequations} \label{eq:reformulation}
\begin{align} \label{eq:dynamics2}
{\rm d}X_t &= f(X_t)\,{\rm d}t + G (I-G^{\rm T} H^{\rm T}C^{-1} H G)^{1/2} {\rm d}\widehat{W}_t + QH^{\rm T} C^{-1/2} {\rm d}\widehat{V}_t,\\
{\rm d}Y_t &= Hf(X_t)\,{\rm d}t + C^{1/2} {\rm d}\widehat{V}_t, \label{eq:obs2}
\end{align}
\end{subequations}
where $\widehat{W}_t$ and $\widehat{V}_t$ denote mutually independent standard Brownian motions of dimension $N_w$ and $N_y$, 
respectively. These equations  correspond exactly to the correlated noise example from \cite[Section 3.8]{sr:BC09}. Furthermore, $H=I$ and $R=0$ lead to
$E_t^{\rm c} = 0$, $QH^{\rm T}C^{-1/2} =C^{1/2}$, and, hence, ${\rm d}X_t = {\rm d}Y_t$.

A straightforward application of the results from \cite[Section 3.8]{sr:BC09} yields the following statement:

\begin{Lemma}[Generalised Kushner--Stratonovich equation] 
\label{lem:KS}	
	The conditional expectations $\pi_t[\phi] = \mathbb{E}[\phi(X_t) \vert Y_{[0,t]}]$ satisfy
\begin{equation}
	\label{eq:KS}
	\begin{aligned}
\pi_t [\phi] & = \pi_0[\phi] + \int_0^t  \pi_s [\mathcal{L}\phi] \,
 \mathrm{d}s  
 + \int_0^t \pi_s \left[\phi h  + HQ\nabla \phi  -\phi \pi_s[h] \right]^{\rm T}
 C^{-1} \left( \mathrm{d}Y_s - \pi_s[h] \, \mathrm{d}s\right), 
	\end{aligned}
\end{equation}
where\footnote{We use the notation $Q:\nabla \nabla \phi = \sum_{i,j=1}^{N_x} Q^{ij} \partial_i \partial_j \phi$.}
\begin{equation}
\label{eq:generator}
\mathcal{L} = f \cdot \nabla + \frac{1}{2} Q : \nabla \nabla
\end{equation}
is the generator of \eqref{eq:dynamics}, $h(x) = Hf(x)$ denotes the observation map, and $\phi$ is a compactly supported smooth function.
\end{Lemma}

For the convenience of the reader, we present an independent derivation in Appendix A. We note that (\ref{eq:KS}) also
arises as the Kushner--Stratonovitch equations for an SDE model (\ref{eq:dynamics}) with observations $Y_t$ satisfying
the observation model
\begin{equation} \label{eq:mod_forward}
{\rm d}Y_t = H \left( f(X_t) - Q \nabla \log \pi_t(X_t)\right) {\rm d}t + C^{1/2}{\rm d}\widetilde{V}_t,
\end{equation}
where $\widetilde{V}_t$ denotes $N_y$-dimensional Brownian motion independent of the Brownian motion $W_t$ in 
(\ref{eq:dynamics}). Here we have used that $\pi_t \left[ HQ \nabla \pi_t \right] = 0$.  This reinterpretation of our state 
estimation problem in terms of uncorrelated model and observation errors and modified observation map
\begin{equation}
\widetilde{h}_t(x) = H \left( f(x) - Q \nabla \log \pi_t(x)\right)
\end{equation}
allows one to apply available MCMC and SMC methods for continuous-time filtering and smoothing problems. See, for example, \cite{sr:stuart15}.
However, there are two major limitations of such an approach. First, it requires approximating the gradient of the log-density. 
Second, the modified observation model (\ref{eq:mod_forward}) is not well-defined in the limit $R\to 0$ and $H=I$, since
the density $\pi_t$ collapses to a Dirac delta function under the given initial condition $X_0 = x_0$ a.s.

In order to circumvent these complications, we develop an alternative approach based on an appropriately modified feedback particle filter 
formulation in the following subsection. 

%
%
\subsection{Generalised feedback particle filter formulation}
%
%

While it is clearly possible to apply the standard feedback particle filter formulations using (\ref{eq:mod_forward}), the following
alternative formulation avoids the need for approximating the gradient of the log-density.

\begin{Lemma}[Feedback particle filter with correlated innovation] Consider the McKean--Vlasov equation
\begin{equation} \label{eq:fpf2Ss}
{\rm d}\widetilde{X}_t = f(\widetilde{X}_t)\,{\rm d}t + G\,{\rm d} W_t + K_t(\widetilde{X}_t) \circ {\rm d} I_t + \Omega_t(\widetilde{X}_t)
\, {\rm d}t,
\end{equation}
where the gain $K_t \in \mathbb{R}^{N_x \times N_y}$ solves 
\begin{equation}
\label{eq:gain equation}
\nabla \cdot \left( \widetilde{\pi}_t \left( K_t C - QH^{\rm T} \right)  \right) = - \widetilde{\pi}_t \left(  h  - \widetilde{\pi}_t[h] 
\right)^{\rm T}, \quad \widetilde{\pi}_t = \Law(\widetilde{X}_t),
\end{equation} 
with observation map $h(x) = Hf(x)$. The function $\Omega_t$ is given by
\begin{equation} \label{eq:def_Omega}
\Omega_t^i = -\frac{1}{2} \sum_{l=1}^{N_x}\sum_{j=1}^{N_y}\partial_l K_t^{ij} (QH^{\rm T})^{lj}, \quad i=1,\ldots,N_x,
\end{equation}
and the innovation process $I_t$ by
\begin{equation} \label{eq:innovation3s}
\mathrm{d}I_t = \mathrm{d}Y_t - \left( h(\widetilde{X}_t) \, \mathrm{d}t + HG \, \mathrm{d} W_t + R^{1/2} \, \mathrm{d}U_t \right).
\end{equation}
Here, $W_t$ and $U_t$ denote mutually independent $N_x$-dimensional and $N_y$-dimensional Brownian
motions, respectively. Then, $\widetilde{\pi}_t = \Law(\widetilde{X}_t)$ coincides with the solution to \eqref{eq:KS}, provided that the initial distributions agree.
\end{Lemma}

It should be stressed that $W_t$ in (\ref{eq:fpf2Ss}) and (\ref{eq:innovation3s}) denote the same Brownian motion, resulting in correlations between the innovation process and model noise.

\begin{proof} In this proof the Einstein summation convention over repeated indices is employed, noting that  \eqref{eq:gain equation} takes the form
	\begin{equation}
	\label{eq:gain equation index}
  \partial_i \left( \widetilde{\pi}_t \left( K_t^{ij}C^{jk} - (QH^{\rm T})^{ik} \right)  \right) = - \widetilde{\pi}_t \left(  h^k  - \widetilde{\pi}_t[h^k] 
	\right), \quad k=1,\ldots, N_y.
	\end{equation}
	We begin by writing \eqref{eq:fpf2Ss} in its It\^{o}-form,
	\begin{equation}
	\label{eq:fpf ito}
	\mathrm{d}\widetilde{X}_t = f(\widetilde{X}_t) \, \mathrm{d}t + G \, \mathrm{d} W_t + K_t (\widetilde{X}_t)\, \mathrm{d}I_t + 
	\widehat{\Omega}_t (\widetilde{X}_t) \, \mathrm{d}t,
	\end{equation}
	where
	\begin{equation}
	\begin{aligned}
	\widehat{\Omega}_t^i &= \Omega_t^i + \frac{1}{2} \left\{ - \left(\partial_l K^{ij}_t\right) (QH^{\rm T})^{lj} + 2\left(\partial_l K_t^{ij}  \right) K_t^{lk} C^{kj}
	\right\} \\
	&= \left(\partial_l K_t^{ij}  \right) \left\{ K_t^{lk} C^{kj} - (QH^{\rm T})^{lj} \right\}
	\end{aligned}
	\end{equation}
	Here we have used that the covariation between $K_t$ and $I_t$ satisfies
	\begin{equation}
	{\rm d} \left\langle K^{ij},I^j\right\rangle_t = \partial_l K_t^{ij}\left( G^{lk}\, \mathrm{d} \left\langle W^k,I\right\rangle_t + K_t^{lk} \,{\rm d} 
	\left\langle I^k,I^j \right\rangle_t \right),
	\end{equation}
	and furthermore $\langle GW,I\rangle_t = -QH^{\rm T} t$ as well as $\langle I, I \rangle_t = 2C t$. 

	For a smooth compactly supported test function $\phi$, It\^o's formula implies
	\begin{equation}
	\label{eq:ito formula}
	\phi(\widetilde{X}_t) = \phi(\widetilde{X}_0) + \int_0^t \partial_i \phi(\widetilde{X}_s) \, \mathrm{d}\widetilde{X}^i_s + \frac{1}{2} \int_0^t \partial_i \partial_j \phi(\widetilde{X}_s) \, \mathrm{d} \langle \widetilde{X}^i, \widetilde{X}^j\rangle_s,
	\end{equation}
	where the covariation process is given by
	\begin{equation}
	\label{eq:covariation}
	\langle \widetilde{X}, \widetilde{X} \rangle_t = tQ - \int_0^t \left( K_sHQ + QH^{\rm T} K_s^{\rm T} \right) \mathrm{d}s +2  \int_0^t K_s C K_s^{\rm T} 
	 \,\mathrm{d}s.
	\end{equation}
	Our aim is to show that $\widetilde{\pi}_t[\phi]$ coincides with $\pi_t[\phi]$ as defined by the Kushner--Stratonovich equation \eqref{eq:KS}. To this end, we insert \eqref{eq:fpf ito} and \eqref{eq:covariation} into \eqref{eq:ito formula} and take the conditional expectation, arriving at
	 \begin{equation}
	 \label{eq:KS from fpf}
	\begin{aligned}
	 \widetilde{\pi}_t [\phi] &=   \widetilde{\pi}_0[\phi] + \int_0^t \widetilde{\pi}_s[\mathcal{L}\phi] \, \mathrm{d}s + \int_0^t \widetilde{\pi}_s\left[(\partial_i \phi) K_s^{ij} \right] \mathrm{d}Y_s^j - \int_0^t \widetilde{\pi}_s\left[(\partial_i \phi) K_s^{ij} h^j \right] \mathrm{d}s \\
	& \qquad + \int_0^t \widetilde{\pi}_s \left[(\partial_i \phi)\, \widehat{\Omega}^i_s \right] \mathrm{d}s +  \int_0^t 
	\widetilde{\pi}_s \left[ \left(\partial_i \partial_j \phi \right)\left( K_s (C K_s^{\rm T} - HQ)  \right)^{ij}\right] \mathrm{d}s,
	 \end{aligned}
	 \end{equation}
	 recalling that the generator $\mathcal{L}$ has been defined in \eqref{eq:generator}. Under the assumption that $K_t$ 
	 satisfies \eqref{eq:gain equation}, the two equations \eqref{eq:KS} and \eqref{eq:KS from fpf} coincide. Indeed, 
	 \begin{equation}
	 \widetilde{\pi}_s\left[(\partial_i \phi) (K_s^{ik}C^{kj} - (QH^{\rm T})^{ij}) \right] 
	 =   \widetilde{\pi}_s \left[ \phi \left( h^j - \widetilde{\pi}_s\left[h^j\right] \right) \right]
	 \end{equation}
	implies
	\begin{equation} \label{eq:Ycontr}
	\widetilde{\pi}_s[\nabla \phi \cdot K_s] = \widetilde{\pi}_s \left[ \phi h + HQ \nabla \phi - \phi \widetilde{\pi}_s[h] \right]^{\rm T} C^{-1},
	\end{equation} 
	and the $\mathrm{d}Y_s$-contributions agree. To verify the same for the $\mathrm{d}s$-contributions,
	we use \eqref{eq:gain equation} to obtain
	\begin{equation}
		\label{eq:calculation1}
		\begin{aligned}
	\widetilde{\pi}_s\left[(\partial_i \phi) K_s^{ij} (h^j-\widetilde{\pi}_t[h^j]) \right] &= - \int_{\mathbb{R}^{N_x}} (\partial_i \phi) K_s^{ij} \partial_l \left( \widetilde{\pi}_s \left( K_s^{ln}C^{nj} - (QH^{\rm T})^{lj}\right)\right)\mathrm{d}x 
	\\
	& = \widetilde{\pi}_s \left[(\partial_i \phi)\, \widehat{\Omega}^i_s \right] +  
	\widetilde{\pi}_s \left[ \left(\partial_i \partial_j \phi \right)\left(K_s (C K_s^{\rm T} - K_sHQ )
	\right)^{ij}\right] .	
	\end{aligned}
	\end{equation}
	Finally, collecting terms in \eqref{eq:KS from fpf} and \eqref{eq:calculation1} and applying (\ref{eq:Ycontr}) to the remaining 
	$\mathrm{d}s$-contribution, i.e. $-\widetilde{\pi}_s [ \nabla \phi \cdot K_s] \widetilde{\pi}_s[h]$, leads to the desired result.
\end{proof}

We note that the correlation between the innovation process $I_t$ and the model error $W_t$ leads to a correction term $\Omega_t$
in (\ref{eq:fpf2Ss}) which cannot be subsumed into a Stratonovitch correction, in contrast to the standard feedback particle filter formulation 
(\ref{eq:fpf1S}). 

\begin{Remark}
Assuming that there exist potential functions $\Psi_t = (\psi_t^1, \ldots, \psi_t^{N_y})$,  $\psi_t^k:\mathbb{R}^{N_x}\to\mathbb{R}$,  solving the Poisson equation(s) \eqref{eq:Poisson} (with $\widetilde{\Pi}_t$ being replaced by $\widetilde{\pi}_t$), 
\eqref{eq:gain equation} can be solved by requiring 
\begin{equation}
K_t = (\nabla \Psi_t + QH^{\rm T}) C^{-1},
\end{equation}
thus generalising \eqref{eq:K Psi}.
\end{Remark}

\begin{Remark} If we set $R=0$, $H=I$, and $K_t = QH^{\rm T}C^{-1} = I$ in (\ref{eq:fpf2Ss}), then one obtains
\begin{equation}
{\rm d}\widetilde{X}_t = {\rm d}Y_t
\end{equation}
since $\Omega_t$ vanishes, and all other terms in (\ref{eq:fpf2Ss}) cancel each other out. If, furthermore, $Y_0 = \widetilde{X}_0 = x_0$ a.s.,
then $\widetilde{X}_t = Y_t$ for all $t \in [0,T]$, which in turn justifies our assumption that the gain $K_t$ is independent of
the state variable. Hence, the McKean--Vlasov formulation (\ref{eq:fpf2Ss}) reproduces the exact reference trajectory $Y_t$ in the case
of no measurement errors and perfectly known initial conditions.
\end{Remark}

We develop a simplified version of the feedback particle filter formulation (\ref{eq:fpf2Ss}) for linear SDEs and Gaussian
distributions in the following subsection, which will form the basis of the generalised 
ensemble Kalman--Bucy filter put forward in the follow-up Section \ref{sec:EnKBF_implementation}.

%
%
\subsection{Generalised Kalman--Bucy filter}
%
%

Let us assume that $f(x) = Fx$ with $F \in \mathbb{R}^{N_x \times N_x}$, that is, equations \eqref{eq:dynamics} 
and \eqref{eq:obs} take the form
\begin{subequations}
	\begin{align}
	\mathrm{d}X_t & = FX_t \,\mathrm{d}t + G \, \mathrm{d} W_t, \\
	\mathrm{d}Y_t & = HF X_t \, \mathrm{d}t + HG \, \mathrm{d}W_t + R^{1/2} \, \mathrm{d}V_t,
	\end{align}
\end{subequations}
with initial conditions drawn from a Gaussian distribution.
In this case $\pi_t$ stays Gaussian for all $t>0$, i.e. $\pi_t \sim {\rm N}(\overline{x}_t, P_t)$ with $\overline{x}_t \in \mathbb{R}^{N_x}$, $P_t \in \mathbb{R}^{N_x \times N_x}$.
Equations \eqref{eq:Poisson}
can be solved uniquely by $\nabla_x \Psi = P_t F^{\rm T} H^{\rm T}$, and thus the McKean--Vlasov 
equations for the feedback particle filter (\ref{eq:fpf2Ss}) reduce to
\begin{equation}
\label{eq:Gaussian feedback}
\mathrm{d}\widetilde{X}_t = F \widetilde{X}_t \, \mathrm{d}t + G \, \mathrm{d} W_t + 
\left( P_t F^{\rm T} H^{\rm T} + QH^{\rm T} \right) C^{-1} \mathrm{d}I_t,
\end{equation}
with the innovation process (\ref{eq:innovation3s}) leading to
\begin{equation}
\label{eq:Gaussian innovation}
\mathrm{d}I_t = \mathrm{d}Y_t - HF\widetilde{X}_t \mathrm{d}t - HG\,\mathrm{d} W_t  - R^{1/2}\mathrm{d}U_t .
\end{equation}
We take the expectation in \eqref{eq:Gaussian feedback}--\eqref{eq:Gaussian innovation} and end up with
\begin{equation}
\mathrm{d}\overline{x}_t = F \overline{x}_t \, \mathrm{d}t + \left( P_t F^{\rm T}  + Q \right) H^{\rm T} C^{-1} \left(  \mathrm{d}Y_t - HF\overline{x}_t \, \mathrm{d}t \right).
\end{equation}
Defining $u_t := \widetilde{X}_t - \overline{x}_t$, we see that
\begin{equation}
\mathrm{d}u_t = Fu_t \, \mathrm{d}t + G \, \mathrm{d} \widetilde{W}_t - \left( P_t F^{\rm T} + Q\right)H^{\rm T} 
C^{-1} \left(HF u_t  \,\mathrm{d}t + H G\,\mathrm{d} W_t + R^{1/2}\mathrm{d}t\right).
\end{equation}
Next we use 
\begin{equation}
\mathrm{d}\left( u_t u_t^{\rm T} \right) = \mathrm{d}u_t u_t^{\rm T} + u_t \mathrm{d}u_t^{\rm T} + \mathrm{d}\langle u,u^{\rm T} \rangle_t
\end{equation}
and $P_t = \mathbb{E} [u_t u_t^{\rm T}]$ to obtain, after some calculations,
\begin{equation}
\mathrm{d}P_t = (F P_t + P_t F^{\rm T})\, \mathrm{d}t - \left( P_t F^{\rm T} + Q  \right) H^{\rm T} 
C^{-1} H \left( F P_t + Q \right) \mathrm{d}t + Q\, \mathrm{d}t. 
\end{equation}
Hence we have shown that our McKean--Vlasov formulation (\ref{eq:Gaussian feedback}) agrees with the standard Kalman--Bucy
filter equations for the mean and the covariance matrix in the correlated noise case \cite{sr:simon}.

%
%
\subsection{Ensemble Kalman--Bucy filter} \label{sec:EnKBF_implementation}
%
%

The McKean--Vlasov equations (\ref{eq:Gaussian feedback}) for linear systems and Gaussian distributions 
suggest approximating the feedback particle filter formulation (\ref{eq:fpf2Ss}) for nonlinear systems by
\begin{equation}
\mathrm{d}\widetilde{X}_t = f (\widetilde{X}_t) \, \mathrm{d}t + G \, \mathrm{d} W_t + 
\left( P_t^{xh} + QH^{\rm T} \right) C^{-1} \mathrm{d}I_t,
\end{equation}
where the innovation process $I_t$ given by (\ref{eq:innovation3s}) as
before. In other words, we approximate the gain matrix $K_t$ in (\ref{eq:fpf2Ss}) by the state independent term 
$\left( P_t^{xh} + QH^{\rm T} \right) C^{-1}$ 
with the covariance matrix $P_t^{xh}$ defined by
\begin{equation}
P_t^{xh} = \widetilde{\pi}_t\left[ (x-\overline{x}_t)(h(x)-\widetilde{\pi}_t[h])^{\rm T}\right] =
 \widetilde{\pi}_t\left[ x \,(h(x)-\widetilde{\pi}_t[h])^{\rm T}\right]
\end{equation}
where $\widetilde{\pi}_t$ denotes the law of $\widetilde{X}_t$.

We can now generalise the ensemble Kalman--Bucy filter formulation (\ref{eq:EnKBF1a}) for the pure parameter estimation problem to the 
state estimation problem with correlated noise. We assume that $M$ initial state values $\widetilde{X}_0^i$ have been sampled
from an initial distribution $\pi_0$ or, alternatively, $X_0^i = x_0$ for all $i = 1,\ldots,M$ in case the initial condition is known exactly. 
These state values are then propagated under the time-stepping procedure
\begin{equation}  \label{eq:EnKBF1x}
\widetilde{X}_{n+1}^i = \widetilde{X}_n^i + \Delta t f(\widetilde{X}_n^i) + \Delta t^{1/2} G \Theta_n^i +
\left(\widehat{P}_n^{xh} + Q H^{\rm T} \right) \left( C +  \Delta t \widehat{P}_n^{hh} \right)^{-1} \Delta I_n^i
\end{equation}
with $\Theta_n^i \sim {\rm N}(0,I)$, step size $\Delta t>0$, empirical covariance matrices
\begin{subequations}
\begin{align}
\widehat{P}_n^{xh} &= \frac{1}{M-1} \sum_{i=1}^M \widetilde{X}_n^i(h(\widetilde{X}_n^i)-\overline{h}_n^M)^{\rm T}, \qquad \overline{h}_n^M = \frac{1}{M}
\sum_{i=1}^M h(\widetilde{X}_n^i),\\
\widehat{P}_n^{hh} &= \frac{1}{M-1} \sum_{i=1}^M h(\widetilde{X}_n^i)(h(\widetilde{X}_n^i)-\overline{h}_n^M)^{\rm T} ,
\end{align}
\end{subequations}
and innovation increments $\Delta I_n^i$ given by 
\begin{equation}
\Delta I_n^i = \Delta Y_n - \Delta t h(\widetilde{X}_n^i)  - \Delta t^{1/2} HG \Theta_n^i - \Delta t^{1/2} R^{1/2} \Xi_n^i, \qquad
\Xi_n^i \sim {\rm N}(0,I).
\end{equation}

The McKean--Vlasov equations of this section form the basis of the methods proposed for the combined state and parameter estimation
problem to be considered next.

%
%
\section{Combined state and parameter estimation} \label{sec:SPE}
%
%

We now return to the combined state and parameter estimation problem and consider the augmented dynamics
\begin{subequations}
	\label{eq:parameter dynamics}
	\begin{align}
	\mathrm{d}X_t & = f(X_t,A_t)\, \mathrm{d}t + G \, \mathrm{d}W_t, \\
	\mathrm{d}A_t&  = 0,
	\end{align}
\end{subequations}
with observations (\ref{eq:obs}) as before. The initial conditions satisfy $X_0 = x_0$ a.s.~and $A_0 \sim \Pi_0$. 
Let us introduce the extended state-space variable $Z_t = (X_t^{\rm T},A_t^{\rm T})^{\rm T}$. 
In terms of $Z_t$, the equations \eqref{eq:parameter dynamics} and \eqref{eq:obs} take the form
\begin{subequations}
	\begin{align}
\mathrm{d} Z_t = \bar{f}(Z) \, \mathrm{d}t + \bar{G} \, \mathrm{d} W_t, \\
\mathrm{d}Y_t = \bar{H}\, \mathrm{d}Z_t + R^{1/2} \, \mathrm{d}V_t,
\end{align}
\end{subequations}
with 
\begin{equation}
\bar{f}(z) = \begin{pmatrix} f(x,a) \\ 0 \end{pmatrix}, \quad \bar{G} = 
\begin{pmatrix}
G & 0 \\
0 & 0
\end{pmatrix}, \quad
\bar{H} =
\begin{pmatrix}
H & 0
\end{pmatrix}.
\end{equation}
Thus we end up with an augmented state estimation problem of the general structure considered in detail 
in Section \ref{sec:state_estimation} already. Below we provide details on some of the necessary modifications.

%
%
\subsection{Feedback particle filter formulation}
%
%
 The appropriately extended feedback particle filter equation (\ref{eq:fpf2Ss}) leads to
\begin{subequations}
	\begin{align}
	\mathrm{d}\widetilde{X}_t & = f(\widetilde{X},\widetilde{A}_t) \, \mathrm{d}t + G \, \mathrm{d} W_t + (\nabla_x \Psi_t(\widetilde{X}_t,\widetilde{A}_t) + 
	QH^{\rm T}) C^{-1} \circ \mathrm{d}I_t   + \Omega_t(\widetilde{X}_t,\widetilde{A}_t), \\
	\mathrm{d}\widetilde{A}_t & = \nabla_a \Psi_t(\widetilde{X}_t,\widetilde{A}_t) C^{-1} \circ \mathrm{d}I_t ,
	\end{align}
\end{subequations}
where (\ref{eq:innovation3s}) takes the form
\begin{equation}
\label{eq:feedback innovation}
\mathrm{d}I_t = \mathrm{d}Y_t - \left( h(\widetilde{X}_t,\widetilde{A}_t) \,\mathrm{d}t  + HG\,\mathrm{d} W_t + R^{1/2}
\mathrm{d}U_t\right)
\end{equation}
with observation map (\ref{eq:forward}) and the correction $\Omega_t$ is given by (\ref{eq:def_Omega}) with $Q$ replaced by $\bar{Q} = \bar{G}\bar{G}^T$
and $H$ by $\bar{H}$.
In the Poisson equation(s) \eqref{eq:Poisson}, $\widetilde{\Pi}_t$ is replaced by $\widetilde{\pi}_t$ denoting the joint density of $(\widetilde{X}_t,\widetilde{A}_t)$.
We also stress that $\Psi_t$ becomes a function of $x$ and $a$ and we distinguish between gradients with respect to $x$ and $a$ using
the notation $\nabla_x$ and $\nabla_a$, respectively.

Numerical implementations of the extended feedback particle filter are demanding due to the need of solving the Poisson
equation(s) (\ref{eq:Poisson}). Instead we again rely on the ensemble Kalman--Bucy filter approximation, which we 
describe next.

%
\subsection{Ensemble Kalman--Bucy filter}
%

We approximate the joint density $\widetilde{\pi}_t$ of $\widetilde{Z}_t$ by an ensemble of particles 
\begin{equation}
\widetilde{Z}_t^i = \begin{pmatrix} \widetilde{X}_t^i \\ \widetilde{A}_t^i \end{pmatrix},
\end{equation}
that is,
\begin{equation}
\widetilde{\pi}_t \approx \frac{1}{M} \sum_{i=1}^M \delta_{\widetilde{Z}_t^i},
\end{equation}
where $\delta_{z'}$ denotes the Dirac delta function centred at $z'$. The initial ensemble satisfies $X_0^i = x_0$ for all $i=1,\ldots,M$, 
and the initial parameter values $A_0^i$ are independent draws from the prior distribution $\Pi_0$.

At the same time, we make the approximation $\widetilde{Z}_t \sim {\rm N}(\overline{z}_t^M, \widehat{P}^{zz}_t)$ when dealing with the 
Kalman gain of the feedback particle filter. Here the empirical mean $\overline{z}_t^M$ has components
\begin{equation}
\overline{x}_t^M = \frac{1}{M} \sum_{i=1}^M \widetilde{X}^i_t, \quad \overline{a}_t^M = \frac{1}{M} \sum_{i=1}^M \widetilde{A}_t^i,
\end{equation}
and the joint empirical covariance matrix is given by
\begin{equation}
\widehat{P}_t^{zz} = \frac{1}{M-1} \sum_{i=1}^M \widetilde{Z}_t^i (\widetilde{Z}_t - \overline{z}_t^M)^{\rm T} = 
\begin{pmatrix}
\widehat{P}^{xx}_t & \widehat{P}_t^{xa} \\
(\widehat{P}_t^{xa})^{\rm T}& \widehat{P}_t^{aa}
\end{pmatrix}.
\end{equation}
As in Section \ref{sec:EnKBF_implementation}, the solution to  \eqref{eq:Poisson} can be  approximated by
	\begin{equation}
	\nabla_x \Psi_t = P_t^{xh}, \quad \nabla_a \Psi_t =P_t^{ah}, 
	\end{equation}
	where the covariance matrices $P^{xh}_t$ and $P^{ah}_t$ are finally estimated by their empirical counterparts
	\begin{subequations}
	\begin{align}
	\widehat{P}^{xh}_t &= \frac{1}{M-1} \sum_{i=1}^M \widetilde{X}_t^i (h(\widetilde{X}_t^i,\widetilde{A}_t^i) - \overline{h}_t^M)^{\rm T}, \\
	\widehat{P}^{ah}_t &= \frac{1}{M-1} \sum_{i=1}^M \widetilde{A}_t^i (h(\widetilde{X}_t^i,\widetilde{A}_t^i) - \overline{h}_t^M)^{\rm T}, 
	\end{align}
	\end{subequations}
	with $\overline{h}^M_t$ defined by
	\begin{equation}
	\overline{h}^M_t = \frac{1}{M}\sum_{i=1}^M h(\widetilde{X}_t^i,\widetilde{A}_t^i).
	\end{equation}
	
Summing everything up, we obtain the following generalised ensemble Kalman--Bucy filter equations
\begin{subequations} \label{eq:KBF2}
	\begin{align} \label{eq:KBF2a}
	\mathrm{d} \widetilde{X}_t^i & = f(\widetilde{X}_t^i, \widetilde{A}^i_t) \, \mathrm{d}t 
	+ G\, \mathrm{d} W_t^i + (\widehat{P}_t^{xh} + QH^{\rm T} ) C^{-1} \, \mathrm{d}I_t^i, \\
	\mathrm{d}\widetilde{A}_t^i & = \widehat{P}_t^{ah} C^{-1} \, \mathrm{d}I_t^i, \label{eq:KBF2b}
	\end{align}
\end{subequations}
where the innovations are given by
\begin{equation} \label{eq:KBF1d}
\mathrm{d}I_t^i = \mathrm{d}Y_t - \left( h(\widetilde{X}_t^i, \widetilde{A}_t^i) \,\mathrm{d}t  +
HG\,\mathrm{d} W_t^i + R^{1/2} \mathrm{d}U_t^i \right), 
\end{equation} 
and $W_t^i$ and $U_t^i$ denote independent $N_x$-dimensional and $N_y$-dimensional, respectively, Brownian motions for $i=1,\ldots,M$.

The interacting particle equations (\ref{eq:KBF2}) can be time-stepped along the lines discussed in 
Section \ref{sec:EnKBF_implementation} for the pure state estimation
formulation of the ensemble Kalman--Bucy filter.

%
\section{Numerical results} \label{sec:Examples}
%

We now apply the generalised ensemble Kalman--Bucy filter formulation (\ref{eq:KBF2})
with innovation (\ref{eq:KBF1d}) to five different model scenarios. 

\begin{figure}[H]
\begin{center}
\includegraphics[width=0.45\textwidth,trim = 0 0 0 0,clip]{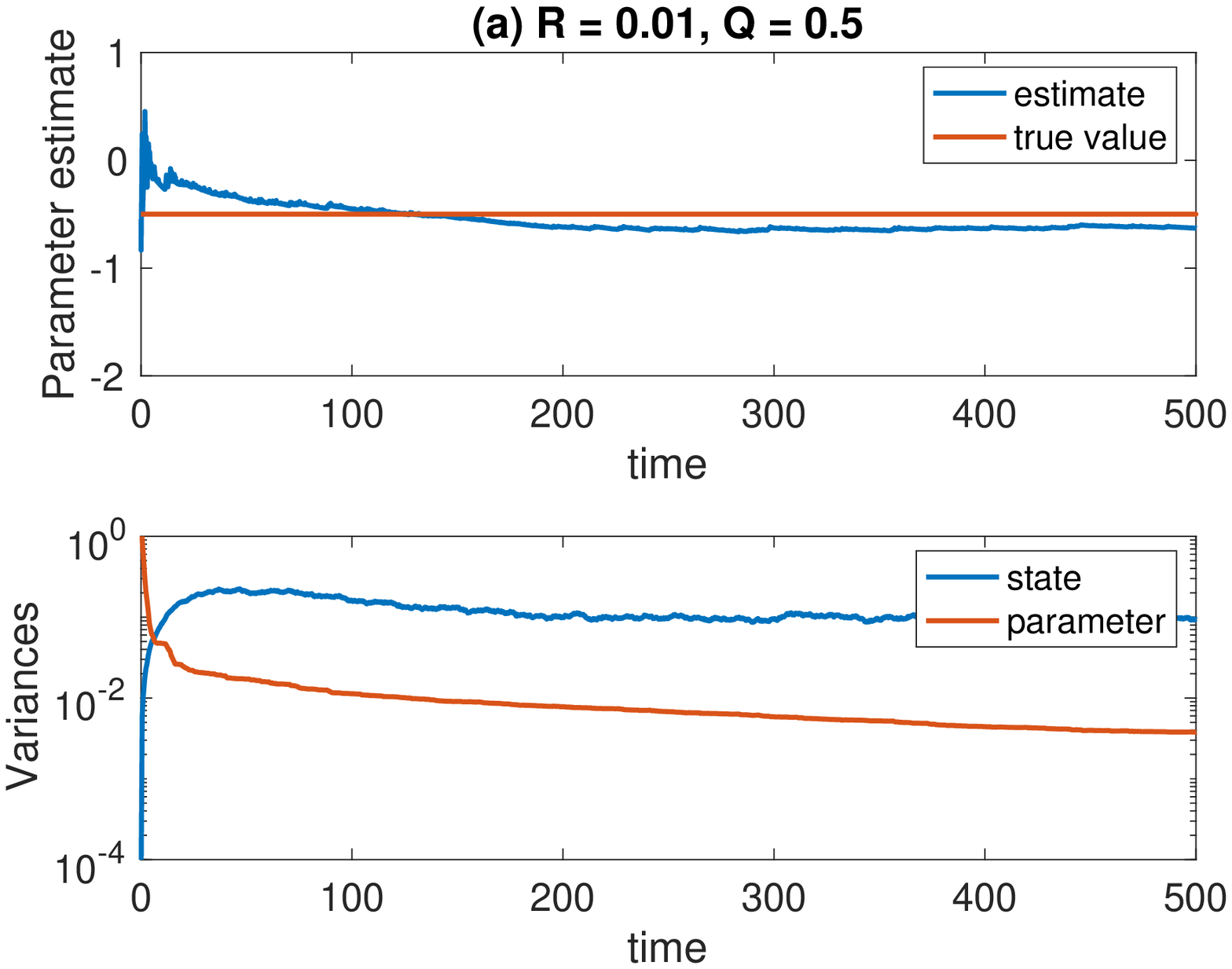} $\qquad$
\includegraphics[width=0.45\textwidth,trim = 0 0 0 0,clip]{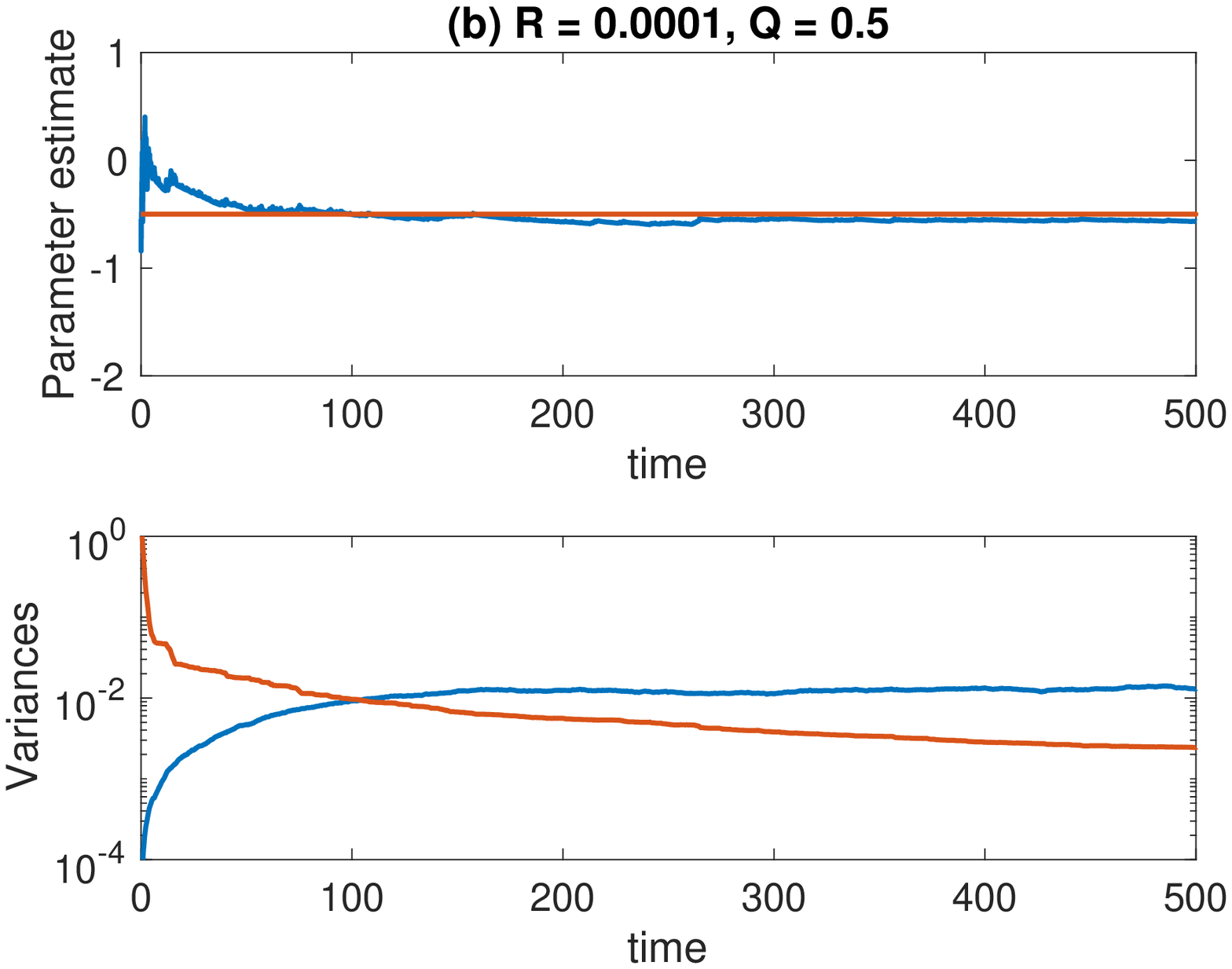} \\ \medskip
\includegraphics[width=0.45\textwidth,trim = 0 0 0 0,clip]{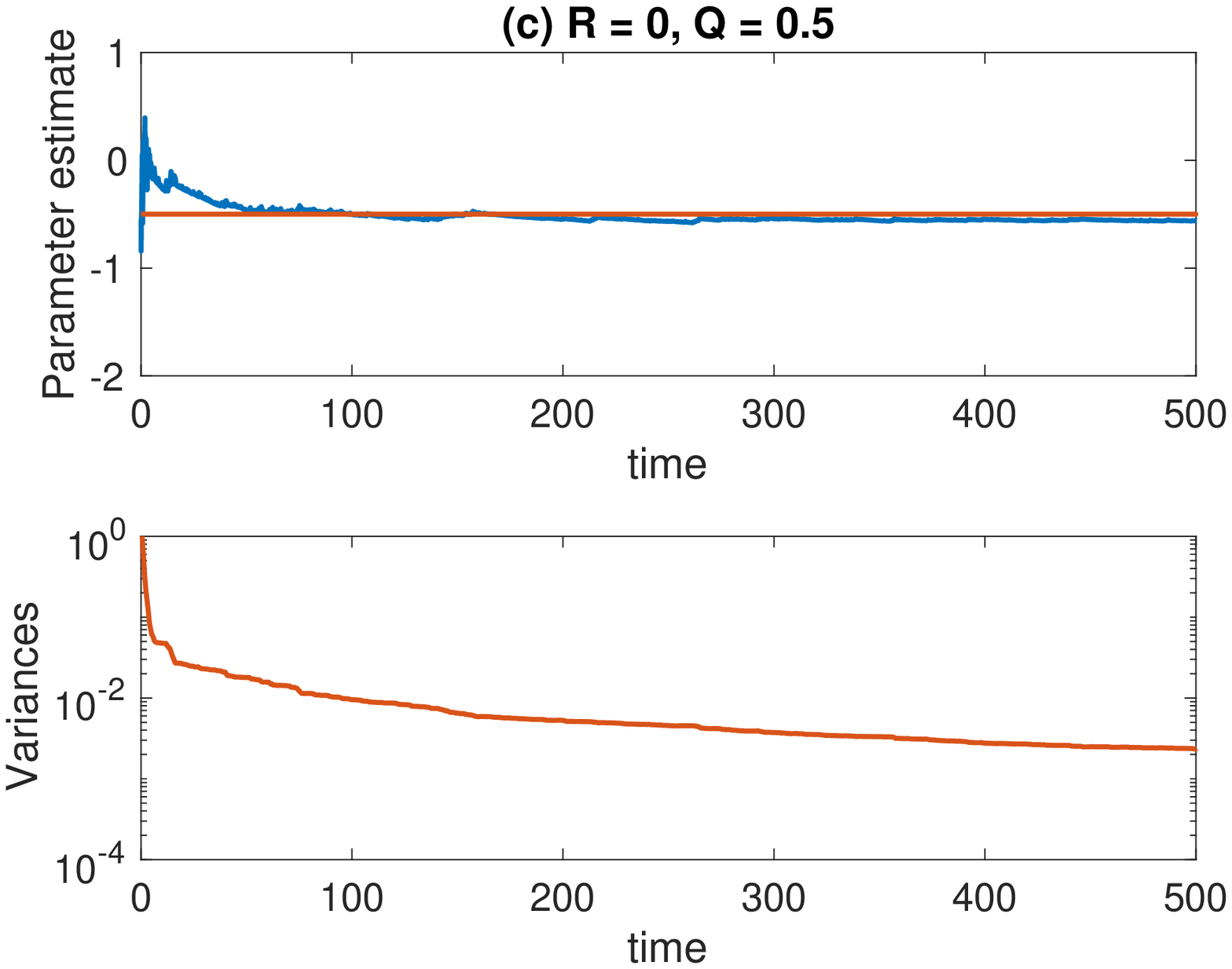} $\qquad$
\includegraphics[width=0.45\textwidth,trim = 0 0 0 0,clip]{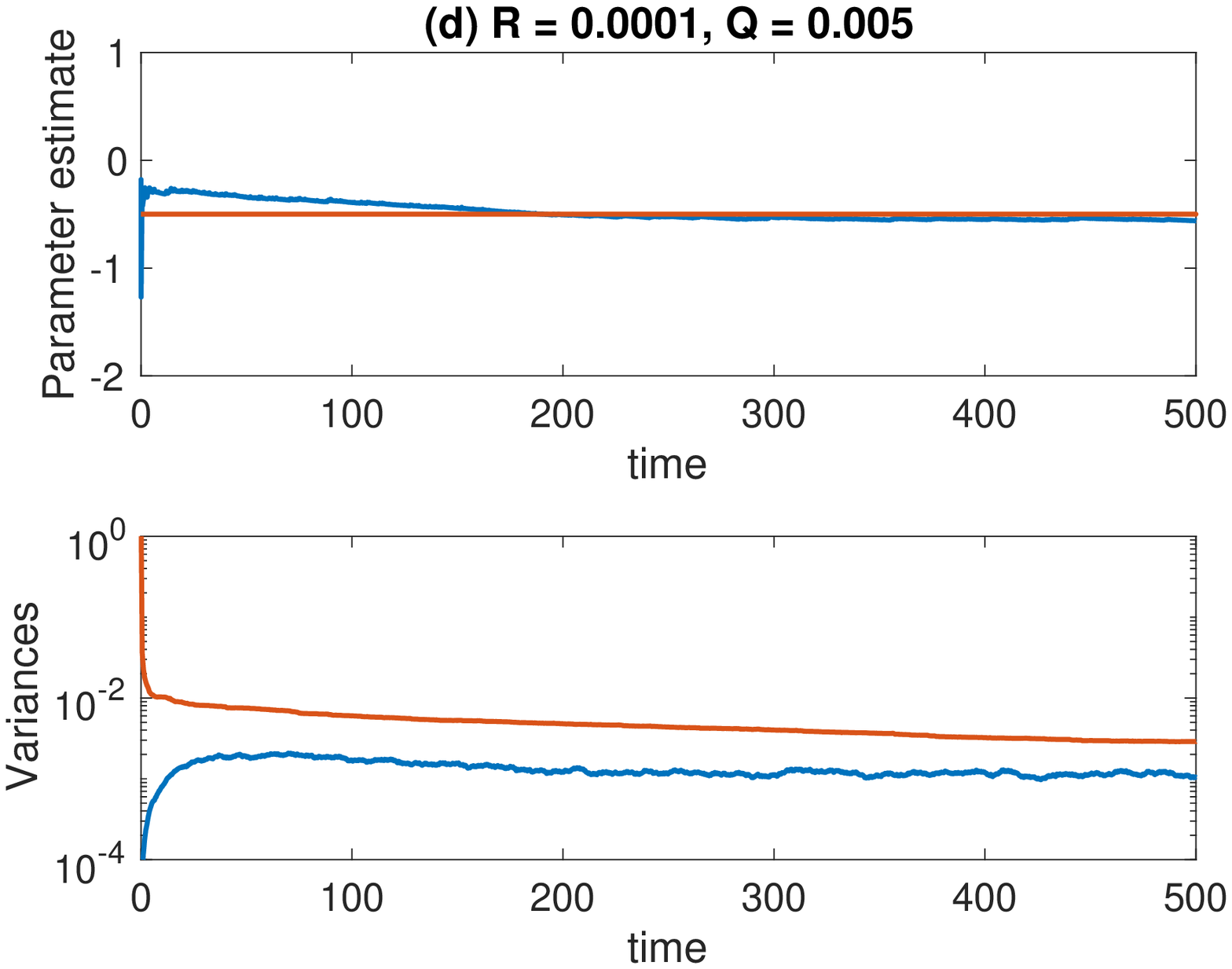} 
\end{center}
\caption{Results for the Ornstein--Uhlenbeck state and parameter estimation problem under different experimental settings:
(a) $Q = 1/2$, $R=0.01$; (b) $Q=1/2$, $R=0.0001$; (c) $Q=1/2$, $R=0$ (pure parameter estimation); (d) $Q=0.005$, $R=0.0001$.
The ensemble size is set to $M=1000$ in all cases. Displayed are the ensemble mean $\overline{a}_n$ and the ensemble variance
in $\widetilde{A}_n$ and $\widetilde{X}_n$. The variance of $\widetilde{X}_n$ is zero when $R=0$ in case (b).}
\label{fig:figure1}
\end{figure}

%
\subsection{Parameter estimation for the Ornstein--Uhlenbeck process}
%

Our first example is provided by the Ornstein--Uhlenbeck process
\begin{equation} \label{eq:example1}
{\rm d}X_t = a X_t\,{\rm d}t + Q^{1/2} {\rm d}W_t
\end{equation}
with unknown parameter $a \in \mathbb{R}$, and known initial condition $X_0 = 1/2$. We assume an observation model of the form 
(\ref{eq:obs}) with $H=1$, and a measurement error taking values $R=0.01$, $R=0.0001$, and $R=0$. The model error variance is set to
either $Q=0.5$ or $Q=0.005$. Except for the case $R=0$ a combined state and parameter estimation problem is
to be solved. We implement the ensemble Kalman--Bucy filter (\ref{eq:KBF2}) with innovation (\ref{eq:KBF1d}),
step size $\Delta t = 0.005$, and ensemble size $M=1000$. The data is generated using the Euler--Maruyama method applied to
(\ref{eq:example1}), with $a = -1/2$ and integrated over a time-interval $[0,500]$ with the same step size.  The prior distribution $\Pi_0$ for the parameter 
is Gaussian with mean $\overline{a} = -1/2$ and variance
$\sigma_a^2 = 2$.  The results can be found in Figure \ref{fig:figure1}. We find that the ensemble Kalman--Bucy filter is able to
successfully identify the unknown parameter under all tested experimental settings, except for the largest measurement error case where $R=0.01$.
There, a small systematic offset of the estimated parameter value can be observed. One can also see that the variance in the parameter
estimate monotonically decreases in time in all cases, while the variance in the state estimates approximately reaches a steady state.

\begin{figure}[H]
\begin{center}
\includegraphics[width=0.45\textwidth,trim = 0 0 0 0,clip]{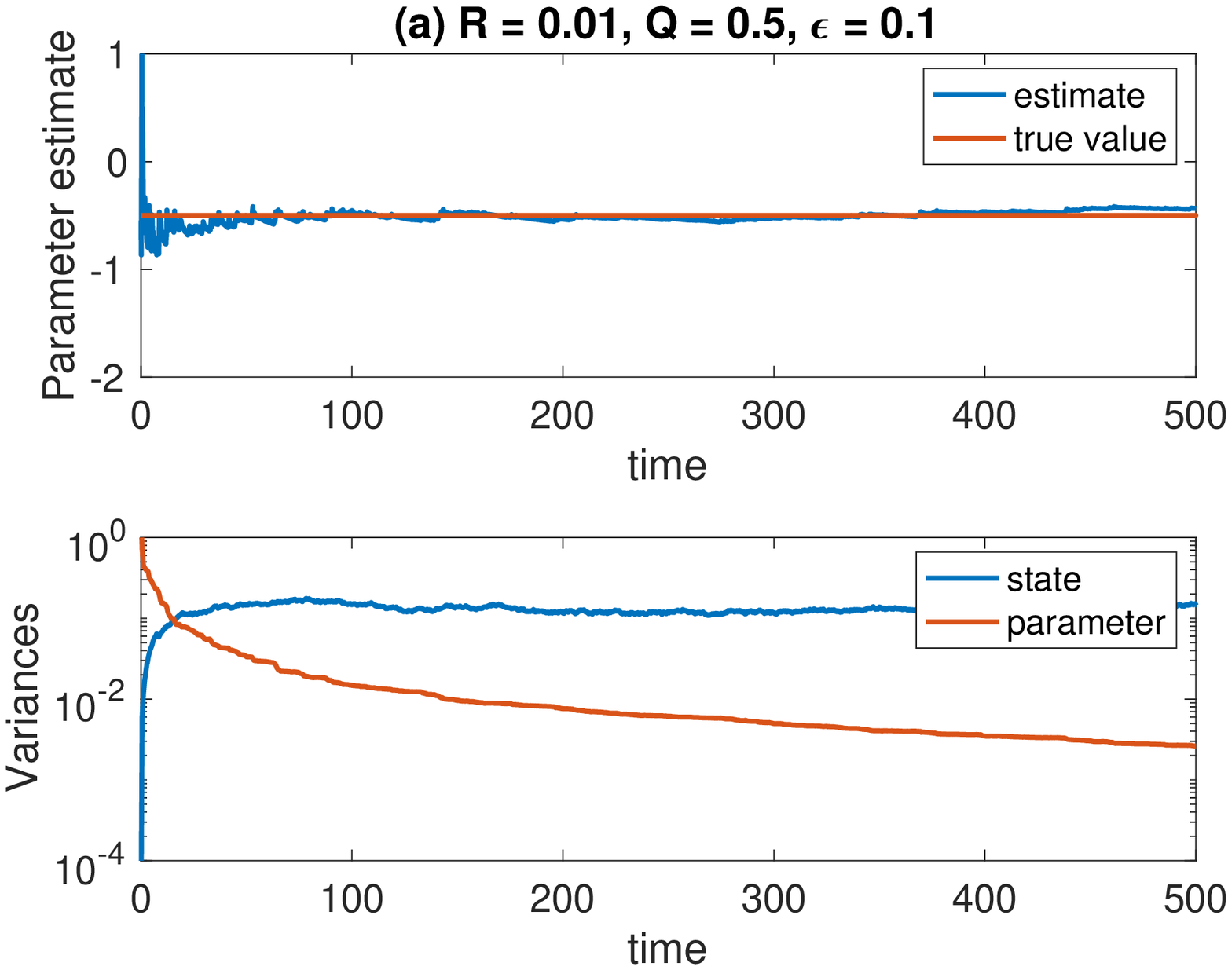} $\qquad$
\includegraphics[width=0.45\textwidth,trim = 0 0 0 0,clip]{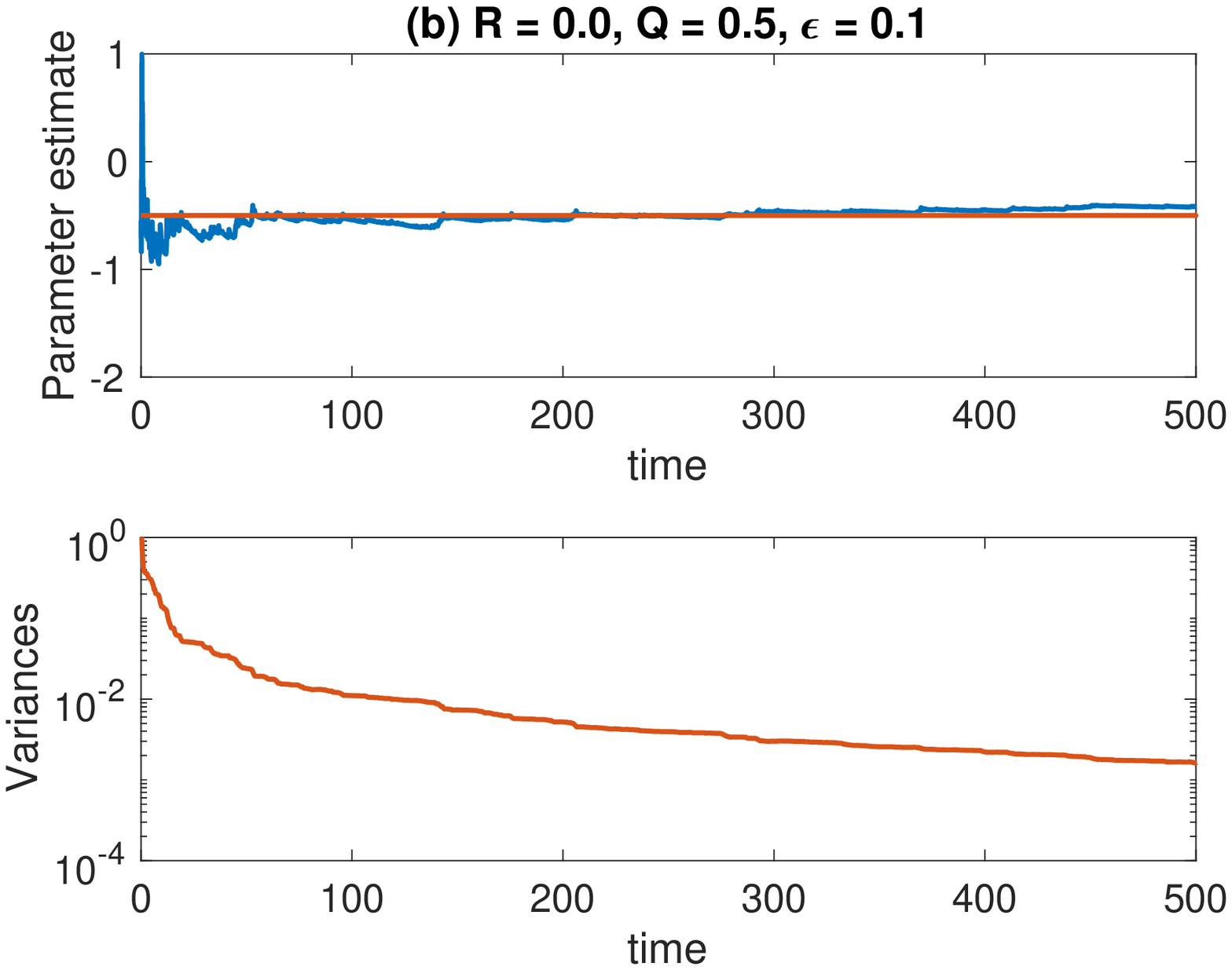} \\ \medskip
\includegraphics[width=0.45\textwidth,trim = 0 0 0 0,clip]{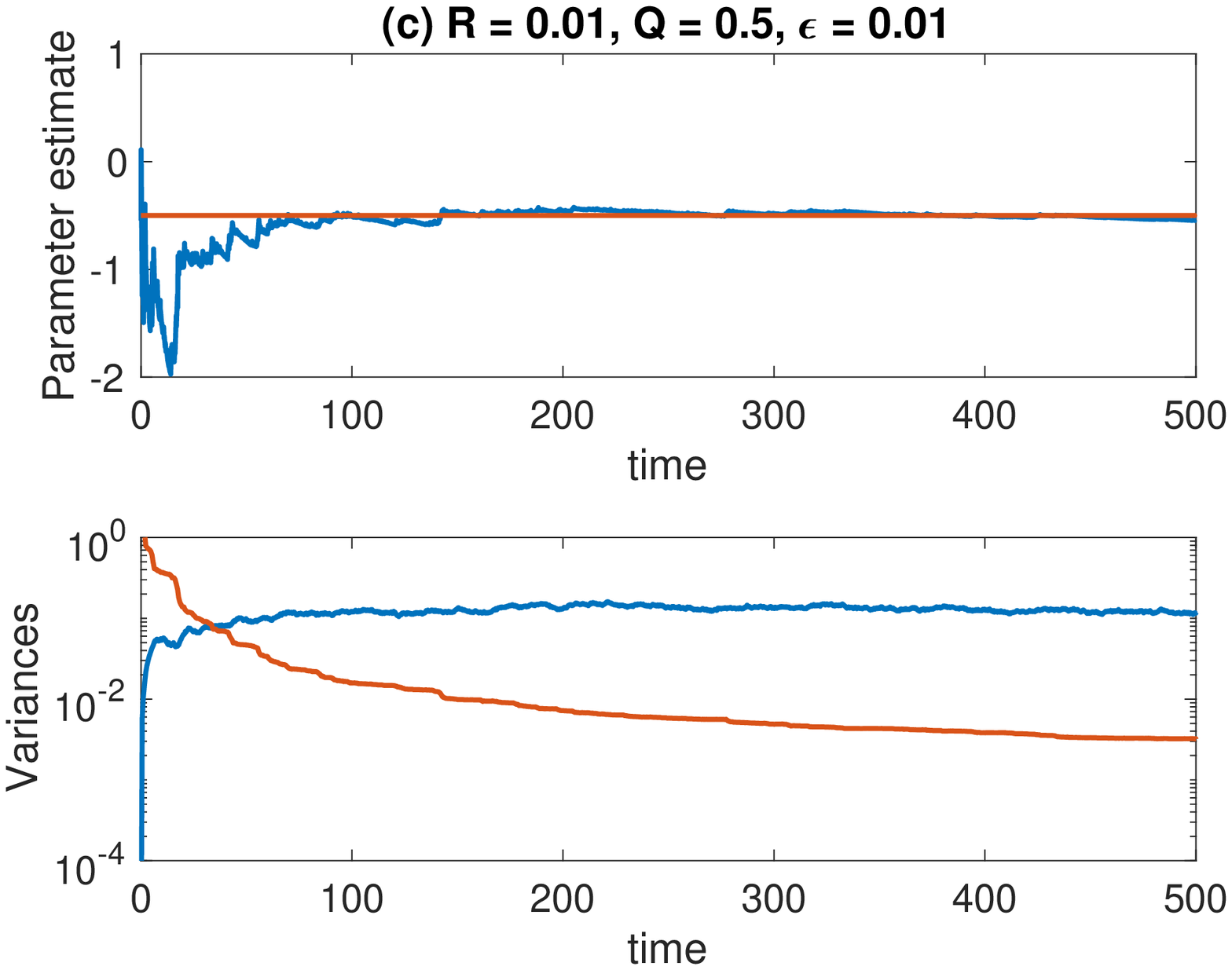} $\qquad$
\includegraphics[width=0.45\textwidth,trim = 0 0 0 0,clip]{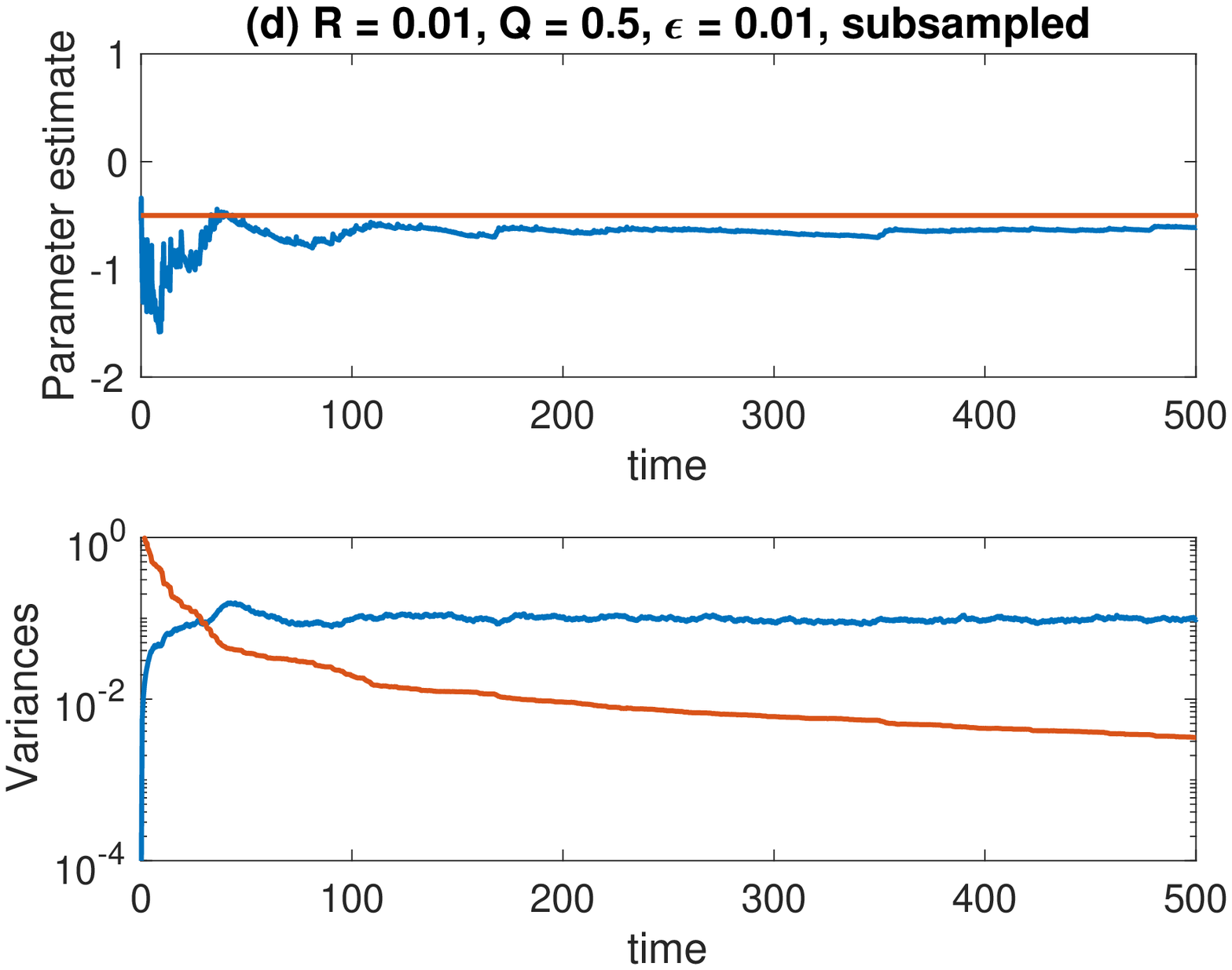} 
\end{center}
\caption{Results for the averaged Ornstein--Uhlenbeck state and parameter estimation problem under different experimental settings:
(a) $Q = 1/2$, $R=0.01$, $\epsilon = 0.1$; (b) $Q=1/2$, $R=0$, $\epsilon = 0.1$ (pure parameter estimation);
(c) $Q=1/2$, $R=0.01$, $\epsilon = 0.01$; (d) $Q=1/2$, $R=0.01$, $\epsilon = 0.01$ and subsampling by a factor of ten.
The ensemble size is set to $M=1000$ in all cases. Displayed are the ensemble mean and the ensemble variance
in $\widetilde{A}_n$ and $\widetilde{X}_n$. The variance of $\widetilde{X}_n$ is zero when $R=0$ in case (b).}
\label{fig:figure2a}
\end{figure}

%
\subsection{Averaging}
%

Consider the equations
\begin{subequations} \label{eq:example2}
\begin{align}
{\rm d}Y_t &= \left( 1-Z_t^2 \right)Y_t \,{\rm d}t + Q^{1/2} {\rm d}W_t^y,\\
{\rm d} Z_t &= -\frac{\alpha}{\epsilon} Z_t \,{\rm d}t + \sqrt{\frac{2\lambda}{\epsilon}} {\rm d}W_t^z 
\end{align}
\end{subequations}
from \cite{sr:ps} for $\lambda,\alpha,\gamma,\epsilon > 0$, and initial condition $Y_0 = 1/2$, $Z_0 = 0$. 
The reduced equations in the limit $\epsilon \to 0$ are given by (\ref{eq:example1}), with parameter value
\begin{equation}  \label{eq:example2b}
a = 1 - \frac{\lambda}{\alpha}
\end{equation}
and initial condition $X_0 = 1/2$. The reduced dynamics corresponds to a (stable) Ornstein--Uhlenbeck 
process for $\lambda/\alpha > 1$. We wish to estimate the parameter $a$ from
observed increments 
\begin{equation} \label{eq:obs_increments}
\Delta Y_n = Y_{n+1} - Y_n + \Delta t^{1/2} R^{1/2} \Xi_n, \qquad \Xi_n \sim {\rm N}(0,1),
\end{equation}
where the sequence of $\{Y_n\}_{n\ge 0}$ is obtained by time-stepping (\ref{eq:example2}) using the
Euler--Maruyama method with a step size $\Delta t$.
We set $\lambda = 3$, $\alpha = 2$ (so that $a = -1/2$), $Q = 0.5$, and $\epsilon \in \{0.1,0.01\}$ in
our experiments. The measurement noise is set to $R = 0.01$ or $R=0$ (pure parameter estimation).

\begin{figure}[H]
\begin{center}
\includegraphics[width=0.45\textwidth,trim = 0 0 0 0,clip]{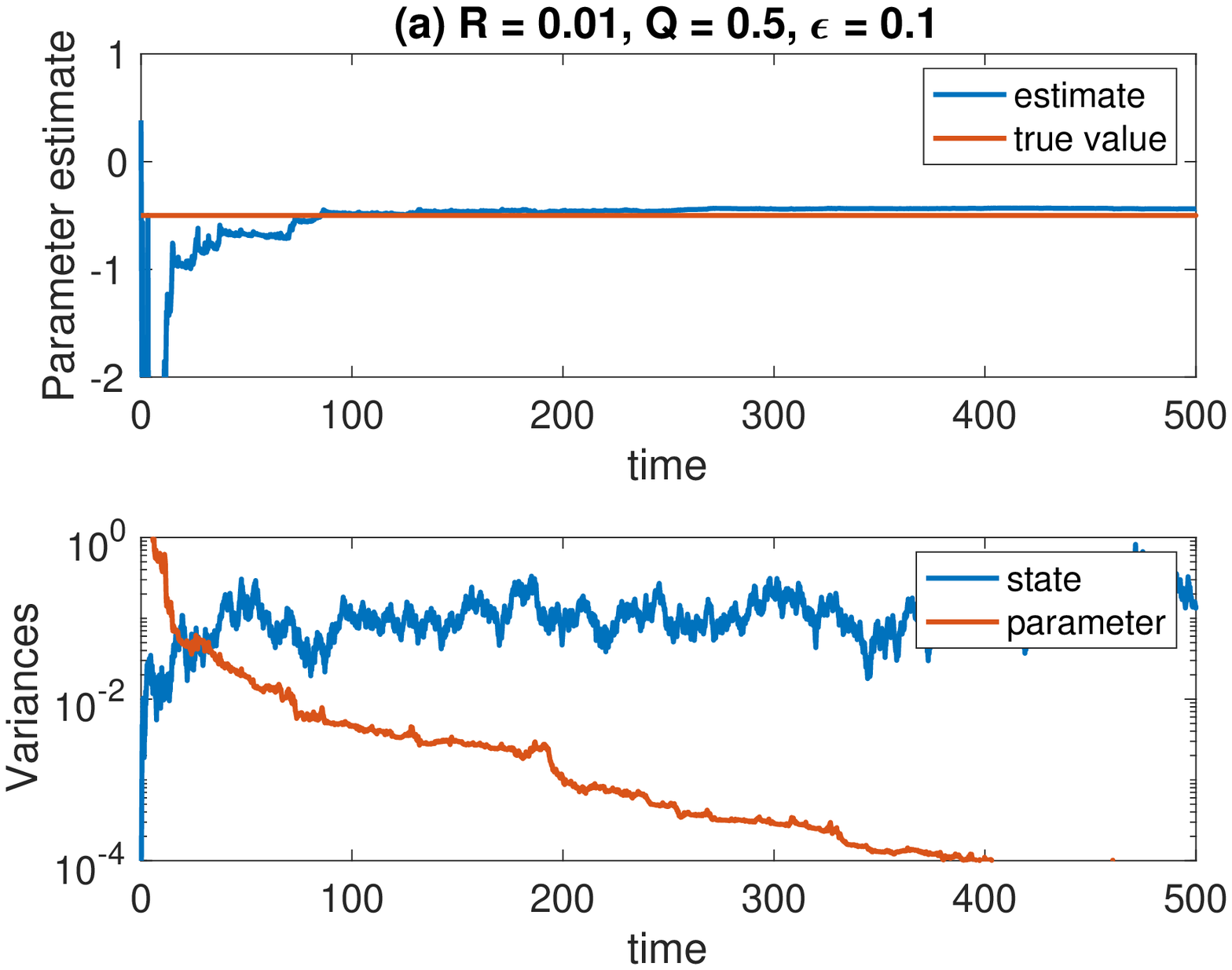} $\qquad$
\includegraphics[width=0.45\textwidth,trim = 0 0 0 0,clip]{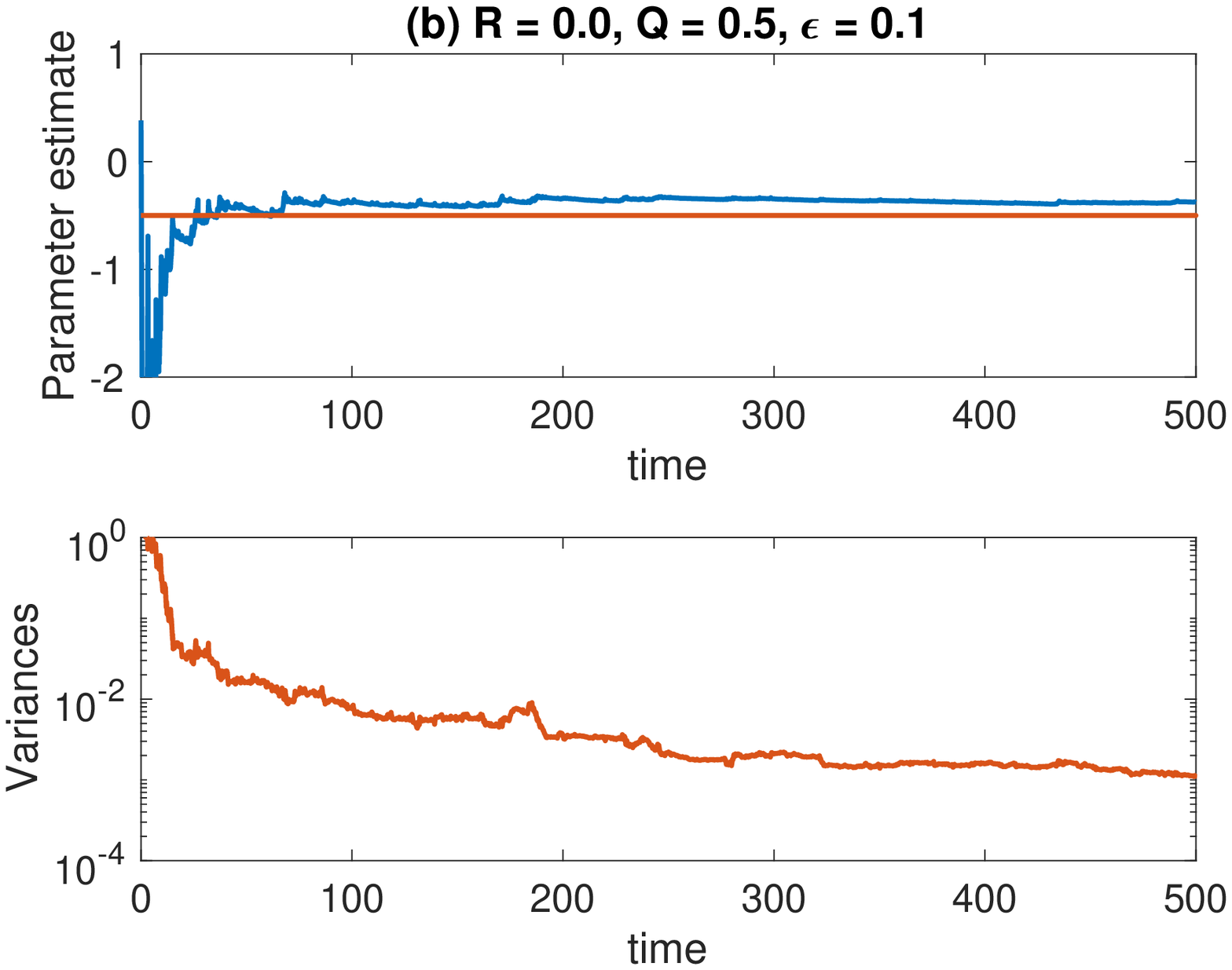} \\ \medskip
\includegraphics[width=0.45\textwidth,trim = 0 0 0 0,clip]{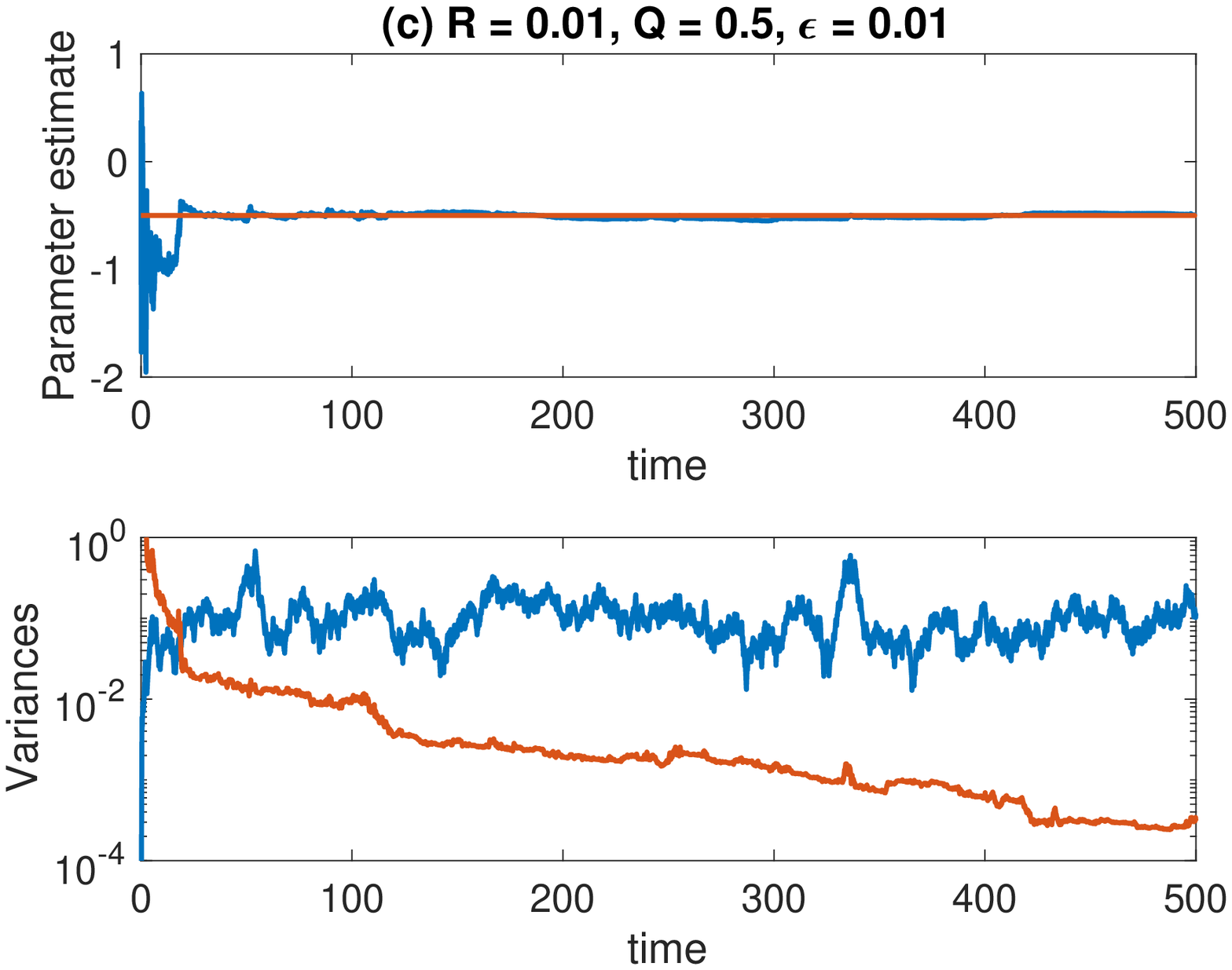} $\qquad$
\includegraphics[width=0.45\textwidth,trim = 0 0 0 0,clip]{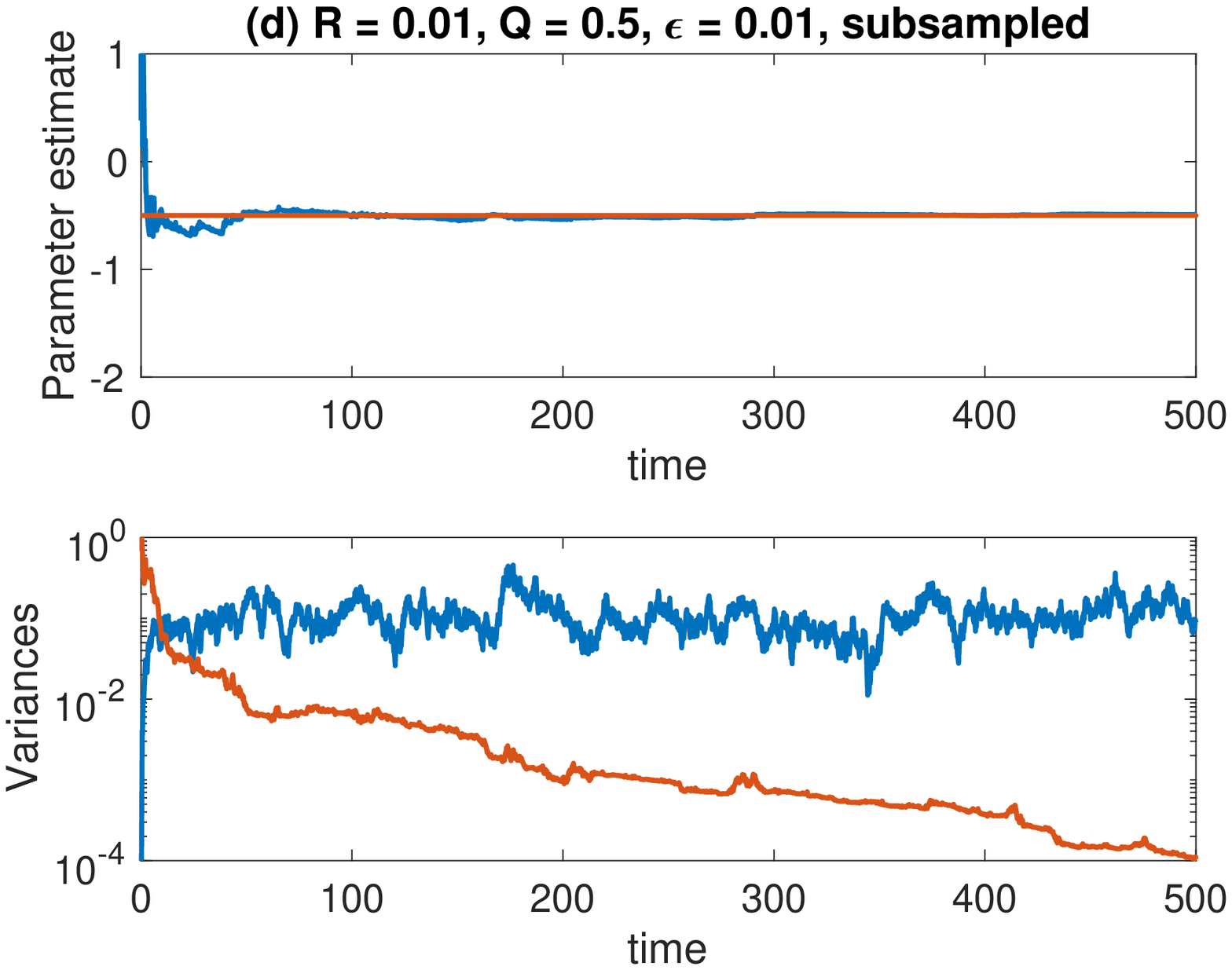} 
\end{center}
\caption{Results as in Figure \ref{fig:figure2a} but with ensemble size $M=10$ instead.}
\label{fig:figure2b}
\end{figure}

We implement the ensemble Kalman--Bucy filter (\ref{eq:KBF2}) with innovation (\ref{eq:KBF1d}),
step size $\Delta t = \epsilon/50$, and ensemble size $M=1000$ for the reduced equations (\ref{eq:example2b}). 
The data is generated from an Euler--Maruyama discretization of (\ref{eq:example2}) with the same step size. 
We also investigate the effect of subsampling the observations for $\epsilon = 0.01$ by solving (\ref{eq:example2}) with
step size $\Delta t = \epsilon/50$ and storing only every tenth solution $Y_n$, while the reduced equations and the ensemble 
Kalman--Bucy filter equations are integrated with $\Delta t = \epsilon/5$. The results are shown in Figure \ref{fig:figure2a}.  Figure \ref{fig:figure2b} shows the results for the same experiments repeated with a smaller ensemble size of $M=10$.
We find that the smaller ensemble size leads to more noisy estimates
for the variance in $\widetilde{X}_n$ and a faster decay of the variance in $\widetilde{A}_n$, 
but the estimated parameter values are equally well converged.   
Subsampling does not lead to significant changes in the estimated parameter values. This is in contrast
to the example considered next.

We finally mention \cite{sr:harlim17} for alternative approaches to sequential estimation in the context of 
averaging using however different assumptions on the data.

\begin{figure}[H]
\begin{center}
\includegraphics[width=0.45\textwidth,trim = 0 0 0 0,clip]{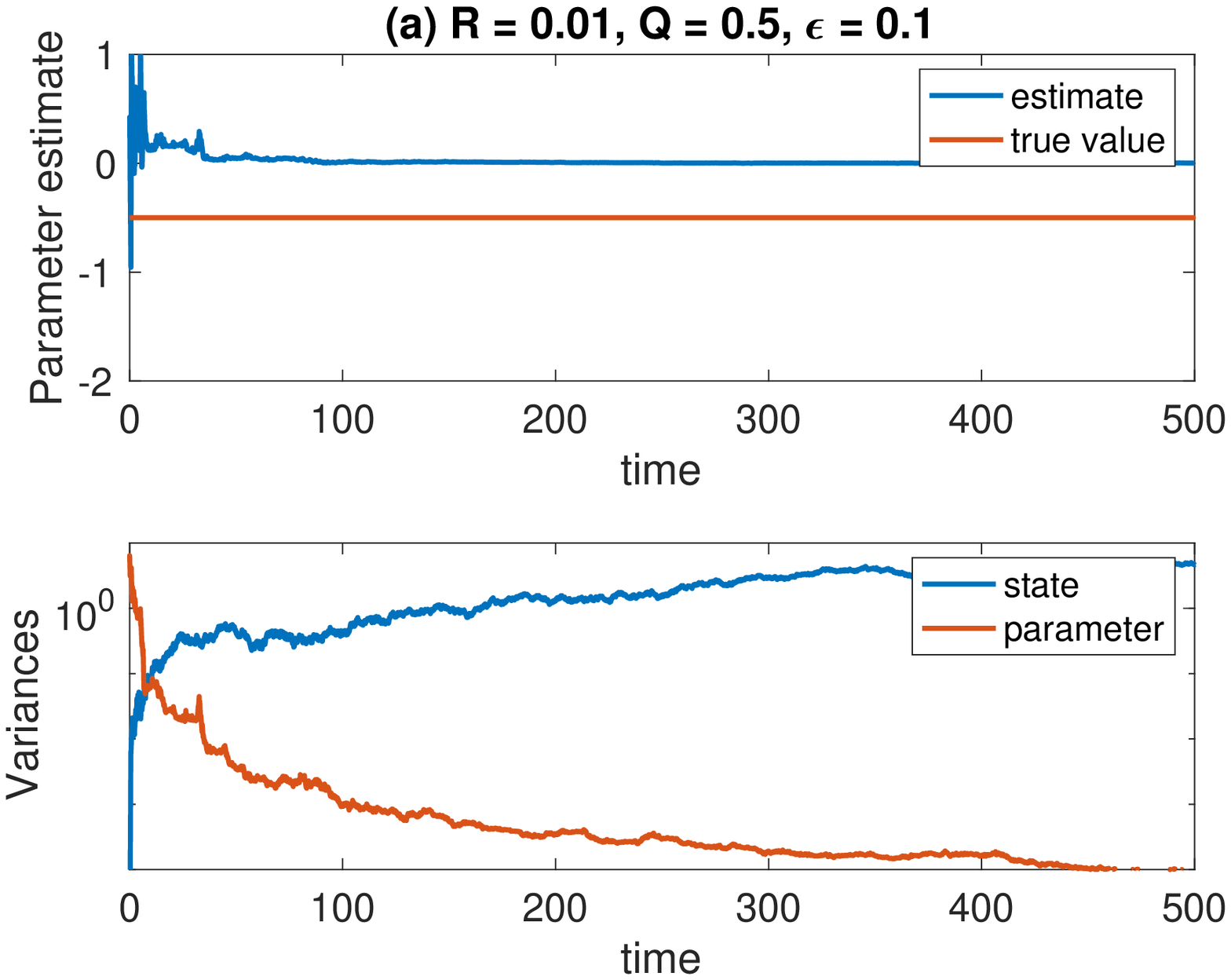} $\qquad$
\includegraphics[width=0.45\textwidth,trim = 0 0 0 0,clip]{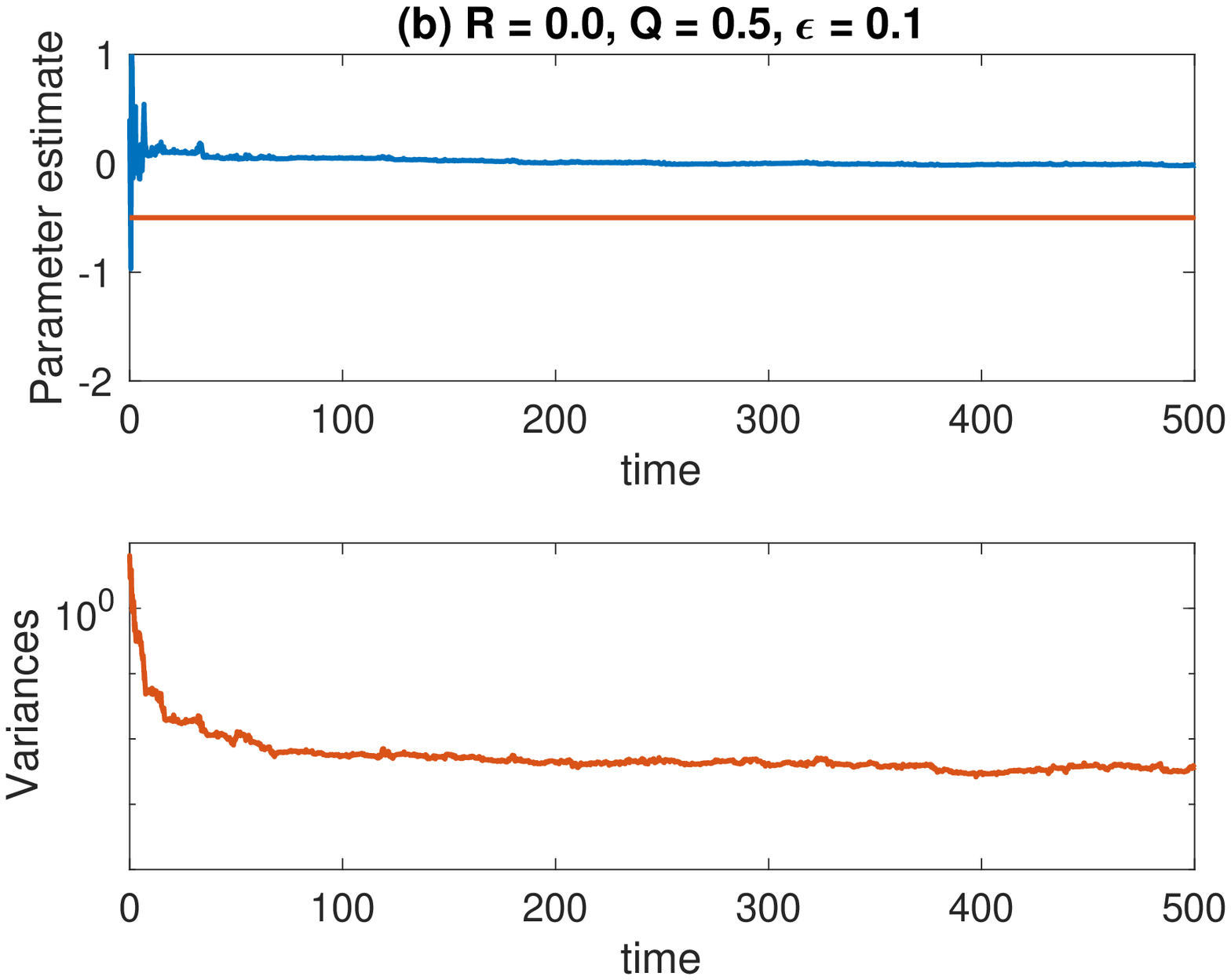} \\ \medskip
\includegraphics[width=0.45\textwidth,trim = 0 0 0 0,clip]{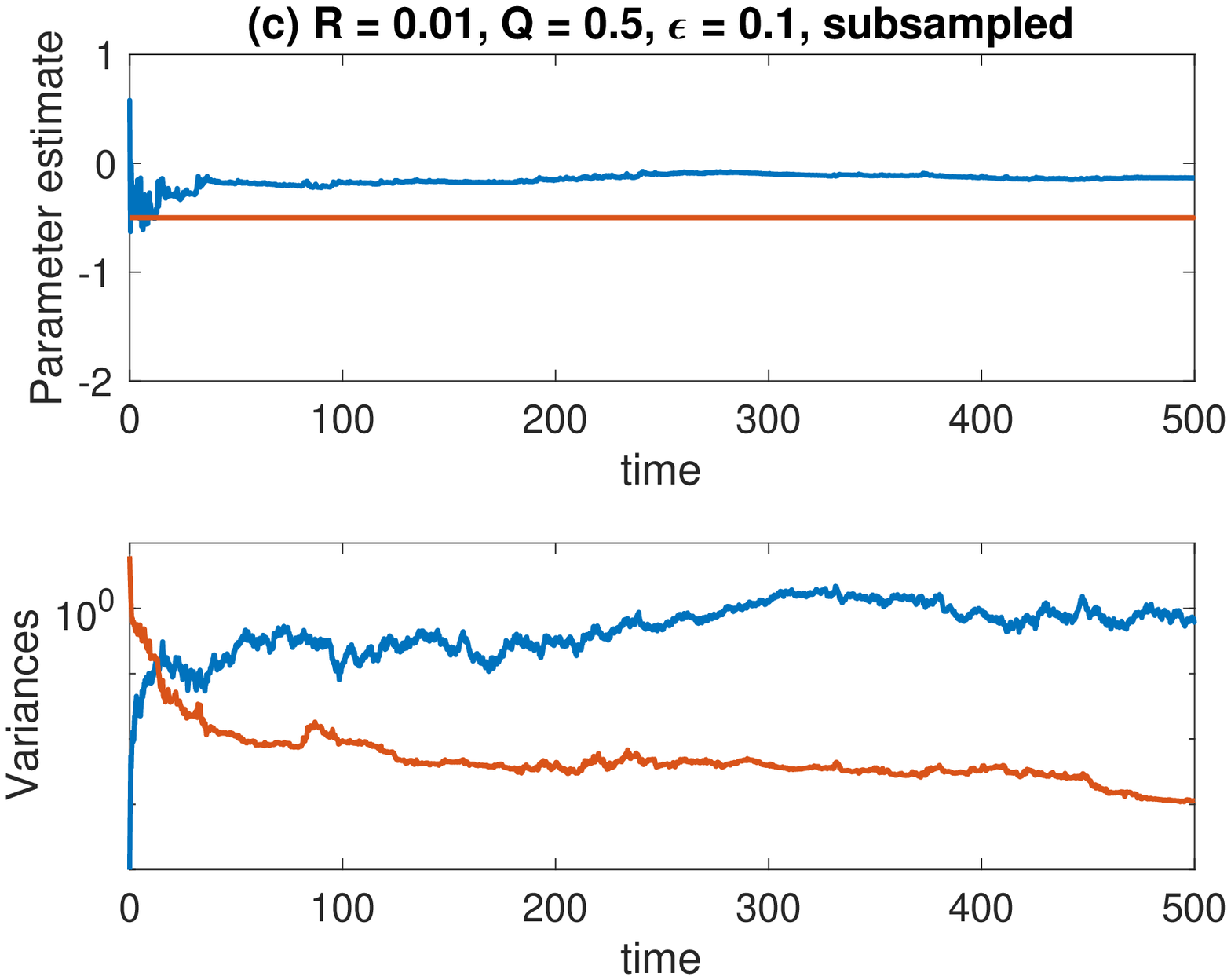} $\qquad$
\includegraphics[width=0.45\textwidth,trim = 0 0 0 0,clip]{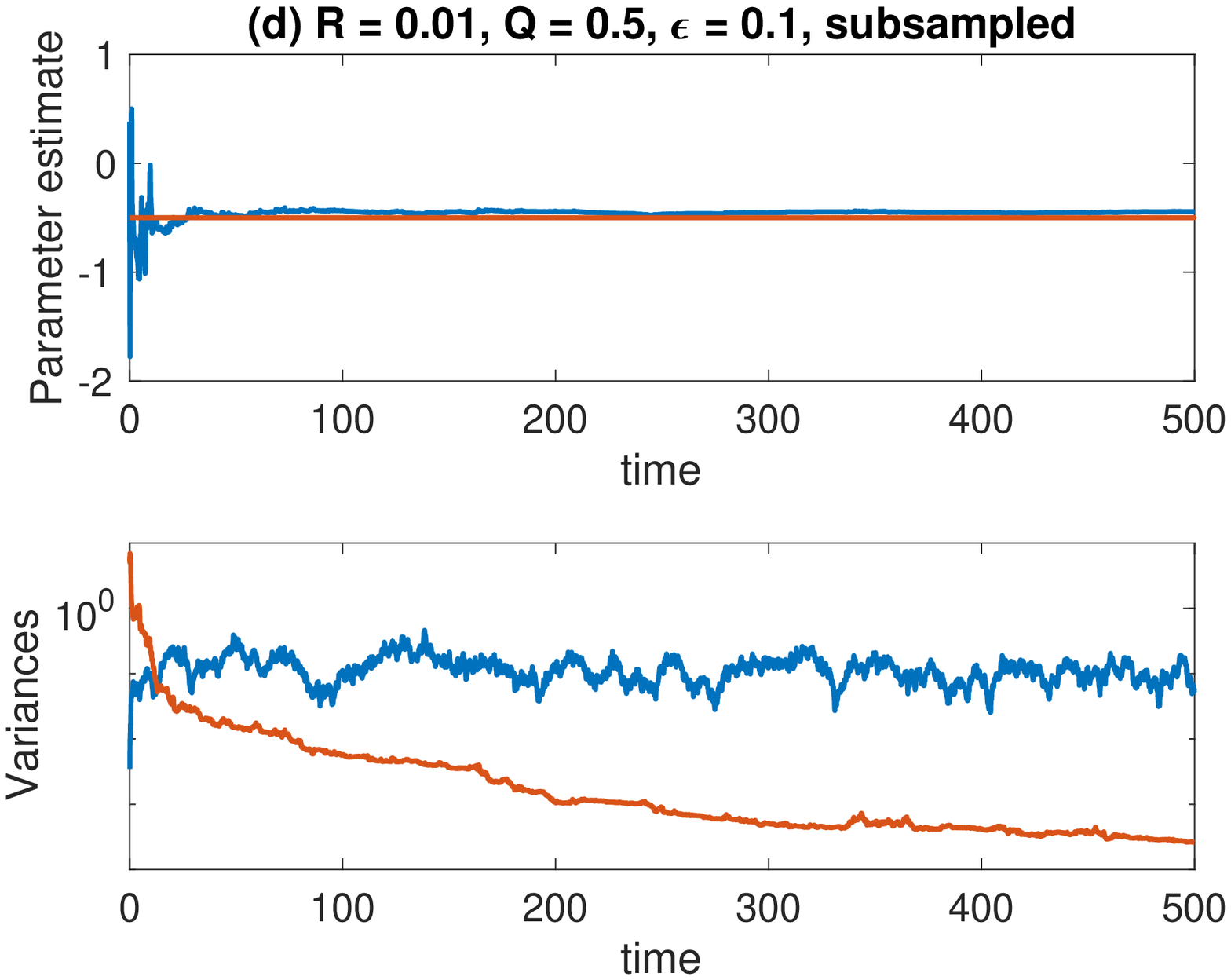} 
\end{center}
\caption{Results for the homoginsation Ornstein--Uhlenbeck state and parameter estimation problem under different experimental settings:
(a) $Q = 1/2$, $R=0.01$, $\epsilon = 0.1$; (b) $Q=1/2$, $R=0$, $\epsilon = 0.1$ (pure parameter estimation);
(c) $Q=1/2$, $R=0.01$, $\epsilon = 0.1$ and subsampling by a factor of fifty; (d) $Q=1/2$, $R=0.01$, $\epsilon = 0.1$ 
and subsampling by a factor of five hundred.
The ensemble size is set to $M=10$ in all cases. Displayed are the ensemble mean and the ensemble variance
in $\widetilde{A}_n$ and $\widetilde{X}_n$. The variance of $\widetilde{X}_n$ is zero under (c).}
\label{fig:figure3}
\end{figure}

%
\subsection{Homogenisation}
%

In this example, the data is produced by integrating the multi-scale SDE
\begin{subequations}
\begin{align} \label{eq:model3a} 
{\rm d}Y_t &= \left( \frac{\sqrt{\sigma/2}}{\epsilon}  Z_t + a Y_t \right) {\rm d} t,\\ \label{eq:model3b}
{\rm d}Z_t &= - \frac{1}{\epsilon^2} Z_t\, {\rm d}t + \frac{\sqrt{2}}{\epsilon} {\rm d} W_t^z
\end{align}
\end{subequations}
with parameter values $\epsilon = 0.1$, $a = -1/2$, $\sigma = 1/2$, and initial condition $Y_0 = 1/2$, $Z_0 = 0$. 
Here, $W_t^z$ denotes standard Brownian motion. The equations are
discretised with step size $\Delta \tau = \epsilon^2/50 = 0.0002$, and the resulting increments (\ref{eq:obs_increments})
are stored over a time interval $[0,500]$. See \cite{sr:KPK11} for more details.

According to homogenisation theory, the reduced model is given by (\ref{eq:example1}) with $Q = \sigma$,
and we wish to estimate the parameter $a$ from the data $\{\Delta Y_n\}$ produced according to (\ref{eq:obs_increments}). 
 It is known that a standard maximum
likelihood estimator (MLE) given by
\begin{equation} \label{eq:ML_estimator}
a_{\rm ML} = \frac{\sum_n Y_{t_n}(Y_{t_{n+1}}-Y_{t_n})}{\sum_n Y_{t_n}^2 \Delta \tau}
\end{equation}
leads to $a_{\rm ML} = 0$ in the limit $\Delta \tau \to 0$ and the observation interval $T\to \infty$. This MLE corresponds to $H=I$ and $R=0$ 
in our extended state space formulation of the problem.  Subsampling can be achieved by choosing an appropriate
time-step $\Delta t > \Delta \tau$ in the ensemble Kalman--Bucy filter equations and a corresponding subsampling of the data
points $Y_n$ in (\ref{eq:obs_increments}). We used $\Delta t = 50\Delta \tau = 0.01$
and $\Delta t = 500 \Delta \tau = 0.1$, respectively. The results can be found in Figure \ref{fig:figure3}. 
It can be seen that only the larger subsampling leads to a correct estimate of the parameter $a$. This is in line with
known results for the maximum likelihood estimator (\ref{eq:ML_estimator}). See \cite{sr:KPK11} and references therein.

\begin{figure}[H]
\begin{center}
\includegraphics[width=0.45\textwidth,trim = 0 0 0 0,clip]{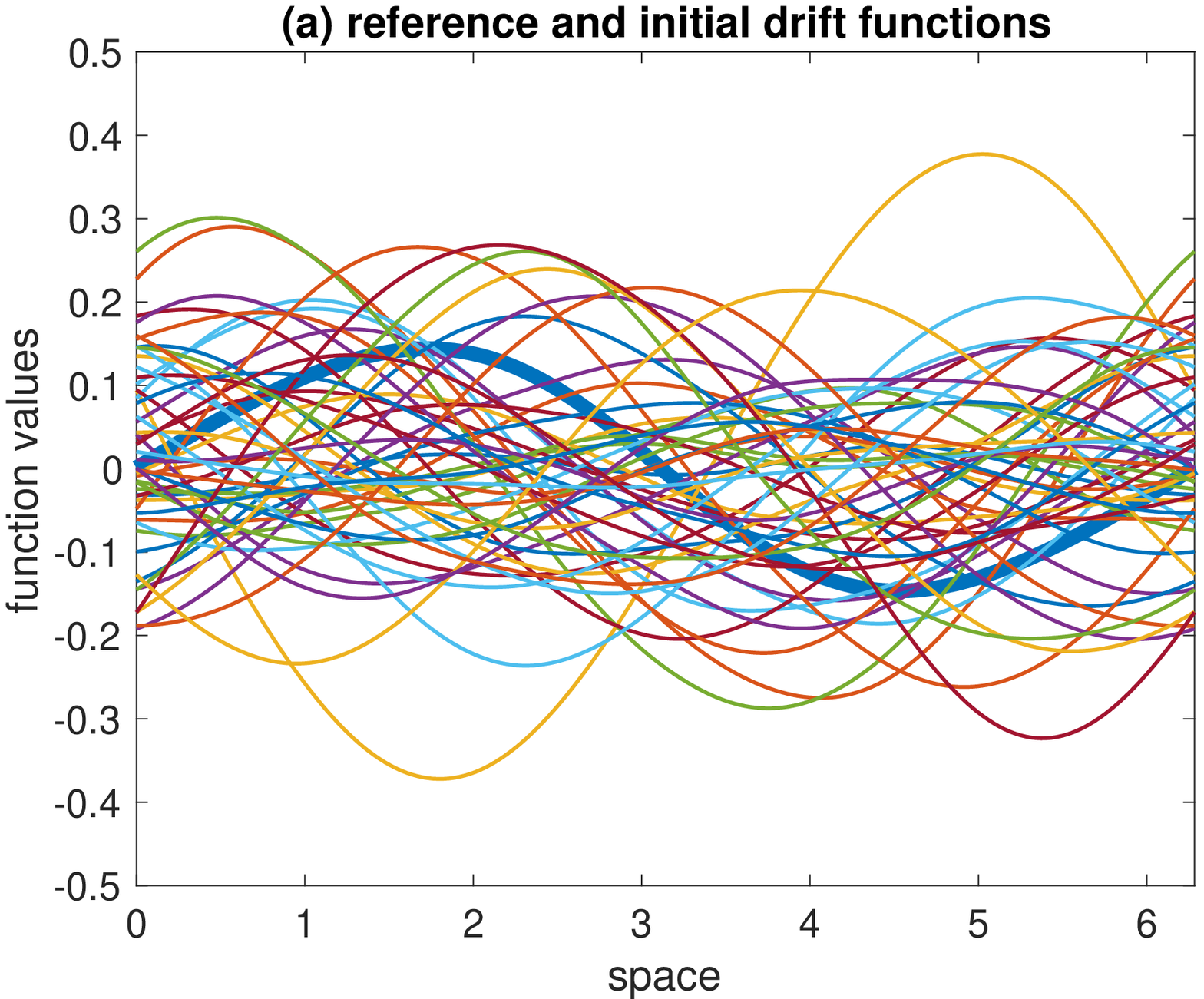} $\qquad$
\includegraphics[width=0.45\textwidth,trim = 0 0 0 0,clip]{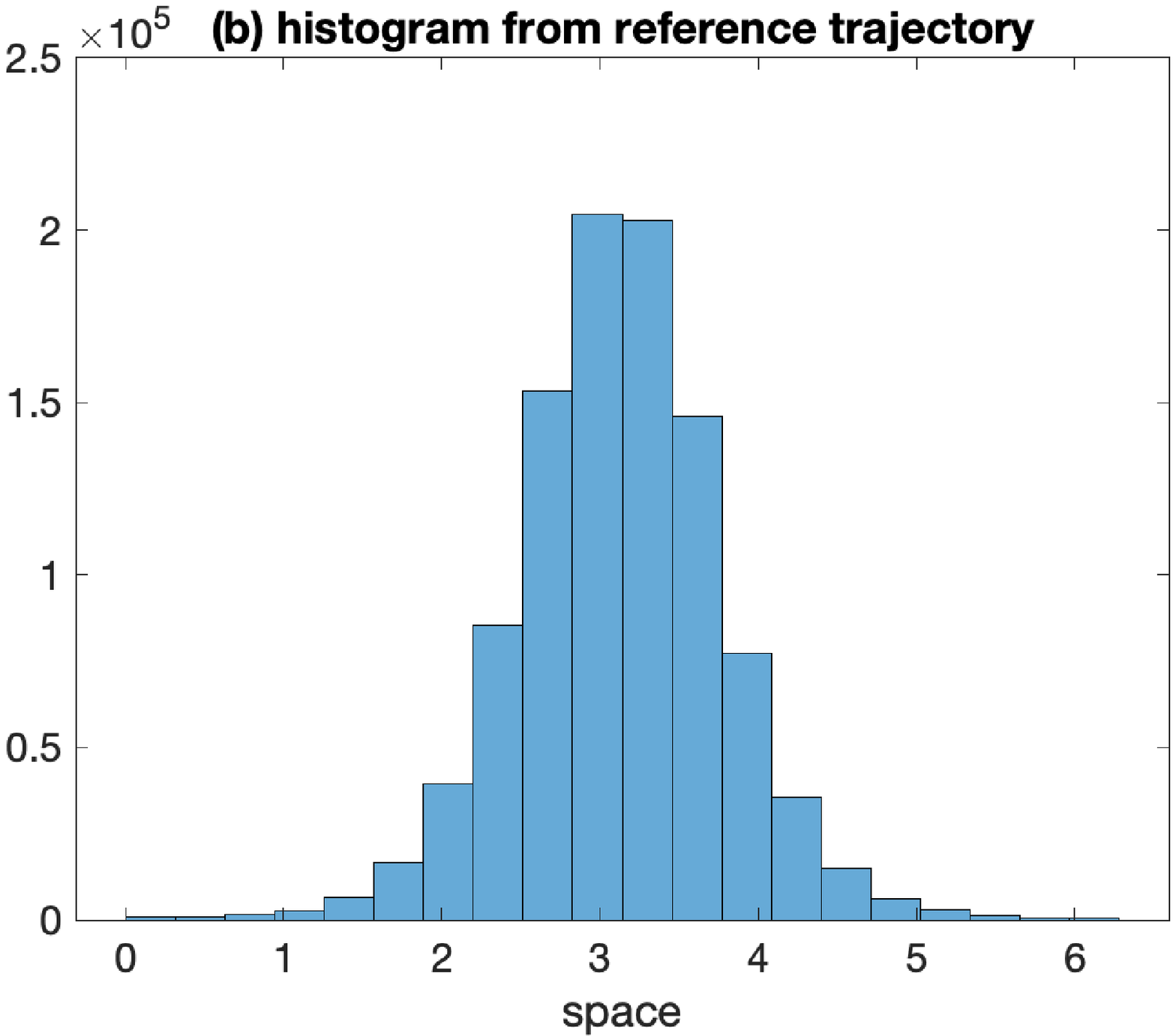} \\ \medskip
\includegraphics[width=0.45\textwidth,trim = 0 0 0 0,clip]{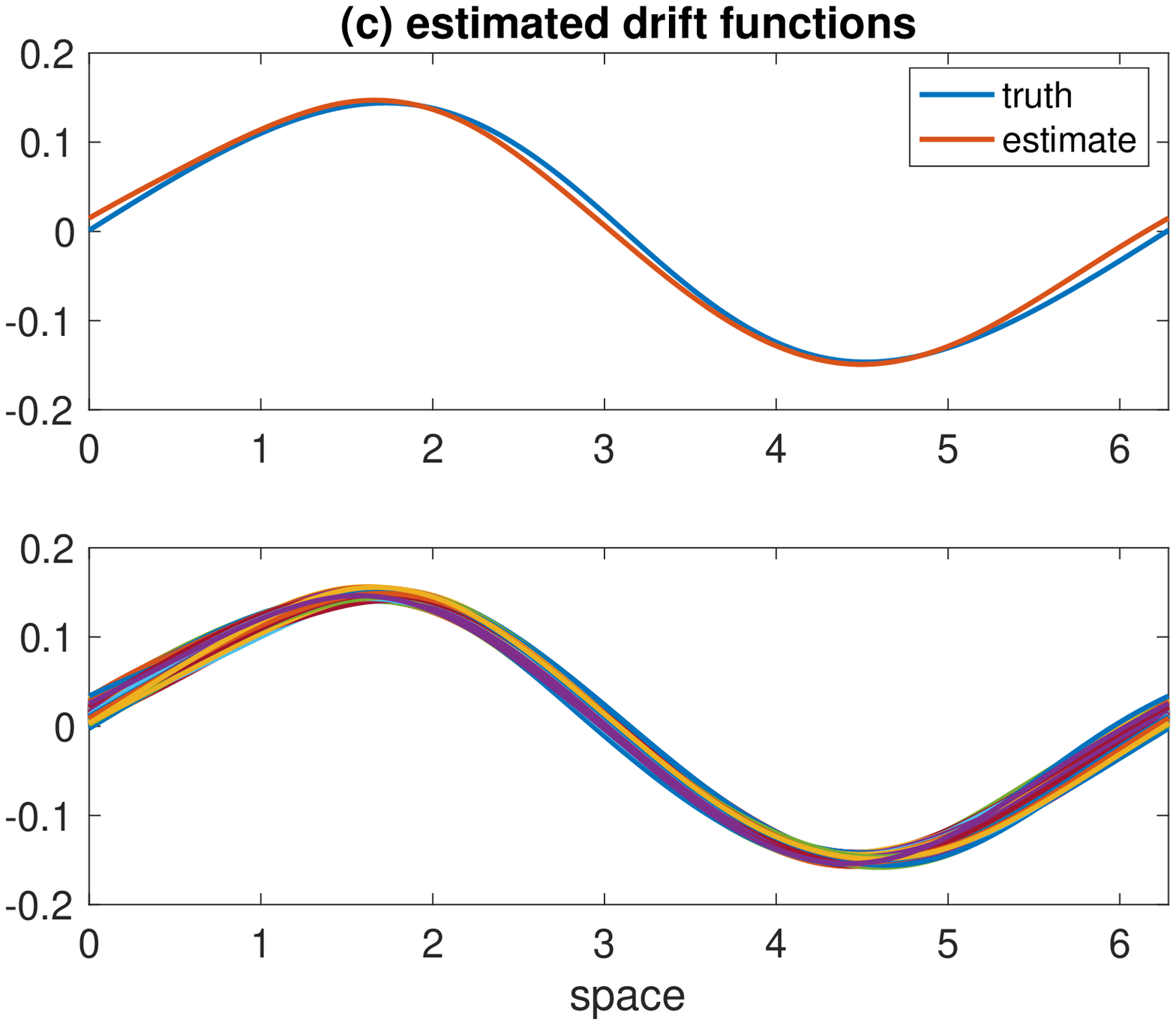} $\qquad$
\includegraphics[width=0.45\textwidth,trim = 0 0 0 0,clip]{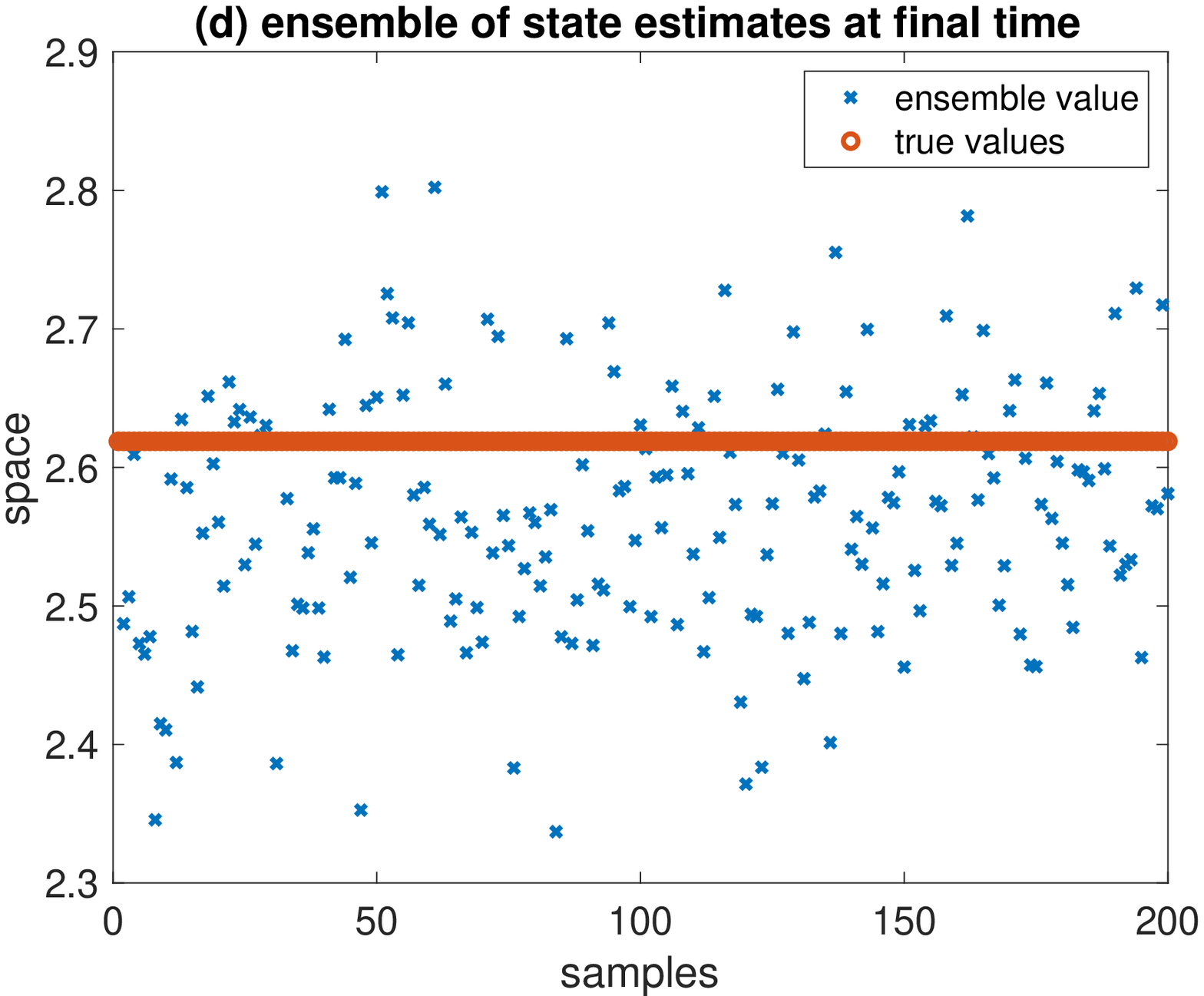} 
\end{center}
\caption{Results for the nonparametric drift and state estimation problem: (a) reference drift function (thick line) 
and ensemble of drift functions drawn from the prior distribution; (b) histogram of samples from the 
reference trajectory; (c) reference drift function and its estimate  (top)
and ensemble of drift functions (bottom) at final time; (d) ensemble of states and the true value at final time.}
\label{fig:figure4}
\end{figure}

%
\subsection{Nonparametric drift and state estimation}
%

We consider nonparametric drift estimation for one-dimensional SDEs over a periodic domain $[0,2\pi)$
in the setting considered from a theoretical perspective in \citep{sr:WaZa16}.  There, a zero-mean Gaussian process prior $\mathcal{GP}(0, \mathcal{D}^{-1})$ is placed on the unknown drift function, with inverse covariance operator
\begin{equation}
\mathcal{D} := \eta [(-\Delta)^{p} + \kappa I].
\end{equation}
The integer parameter $p$ sets the regularity of the process, whereas $\eta, \kappa \in \mathbb{R}^+$ control its characteristic correlation length and stationary variance.

Spatial discretization of the problem is carried out by first defining a grid of $N_d$ evenly spaced points on the domain, at locations $x_i = i\Delta x$, $\Delta x = 2\pi/N_d$.  The drift function is projected onto compactly supported functions centred at these points, which are piecewise linear with
\begin{equation}
  b_i(x_j) = \delta_{ij}
\end{equation}
and linear interpolation is used to define a drift function $f(x,a)$ for all $x \in [0,2\pi)$, that is, it is of the form (\ref{eq:linear_parameter}) with
$f_0(x) \equiv 0$. In this example, we set $N_d = 200$.  Sample realisations, as well as the reference drift $f^*$, 
can be found in Figure \ref{fig:figure4}(a).

Data is generated by integrating the SDE \eqref{eq:dynamics} with drift $f^\ast$ forward in time from initial condition $X_0 = \pi$ and with noise level $Q = 0.1$, using the 
Euler--Maruyama discretisation with step size $\Delta t = 0.1$ over one million time-steps.  The spatial distribution of the solutions $X_n$
is plotted in Figure \ref{fig:figure4}(b).  The data is then given by
\begin{equation}
\Delta Y_n = X_{n+1}-X_n + \Delta t^{1/2} R^{1/2}\Xi_n
\end{equation}
with $R = 0.00001$.  Data assimilation is performed using the time-discretised ensemble Kalman--Bucy filter equations 
(\ref{eq:KBF2}) with innovation (\ref{eq:KBF1d}), ensemble size $M=200$, and step size $\Delta t = 0.1$.

The final estimate of the drift function (ensemble mean) and the ensemble of
drift functions can be found in Figure \ref{fig:figure4}(c). Figure \ref{fig:figure4}(d) displays the ensemble of state estimates and the
value of the reference solution at the final time. We find that the ensemble Kalman--Bucy filter is able to successfully estimate the drift
function and the model states. Further experiments reveal that the drift function can only be identified for sufficiently small measurement errors.

\begin{figure}[H]
\begin{center}
\includegraphics[width=0.45\textwidth,trim = 0 0 0 0,clip]{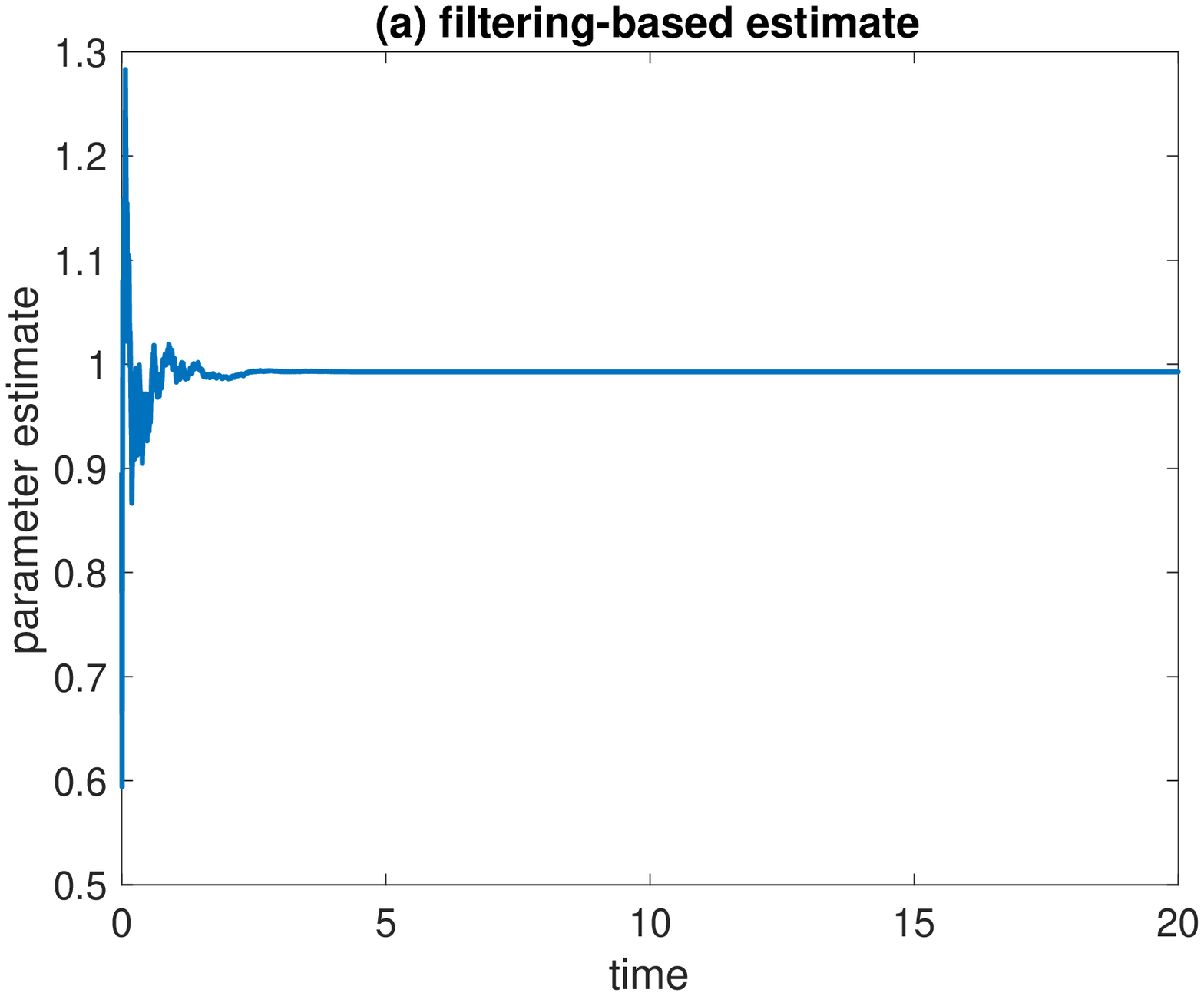} $\qquad$
\includegraphics[width=0.45\textwidth,trim = 0 0 0 0,clip]{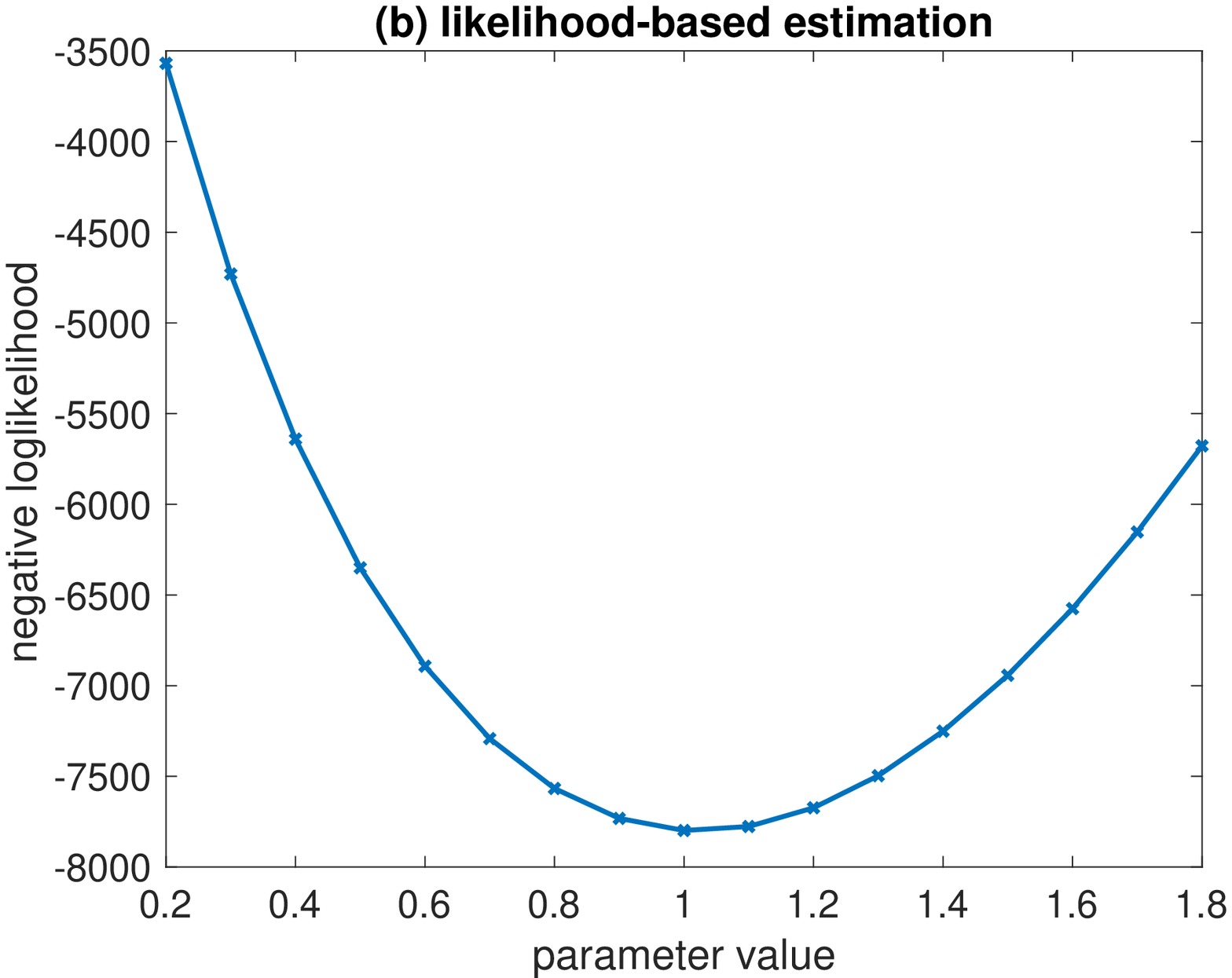} 
\end{center}
\caption{Results for SPDE parameter estimation: (a) estimate of $\theta$ as a function of time as obtained by the
ensemble Kalman--Bucy filter; (b) evidence based on a Kalman--Bucy filter for state estimation 
applied to a sequence of parameter values $\theta \in \{0.2,0.3,\ldots,1.8\}$.}
\label{fig:figure5}
\end{figure}

%
\subsection{SPDE parameter estimation}
%
Consider the stochastic heat equation on the periodic domain $x \in [0,2\pi)$, given in conservative form by the stochastic partial
differential equation (SPDE)
\begin{align}
  \label{eq:spde_fulleqn}
  \dif u(x,t) = \nabla \cdot \left(\theta(x) \nabla u(x,t)\right) \dif t + \sigma^{1/2} \,\dif W(x,t),
\end{align}
where $W(x,t)$ is space-time white noise.  With constant $\theta(x) = \theta$, this SPDE reduces to
\begin{align}
  \label{eq:spde_constthetaeqn}
  \dif u(x,t) = \theta \Delta u(x,t) \dif t + \sigma^{1/2}\, \dif W(x,t).
\end{align}
In this example, we examine the estimation of $\theta$ from incremental measurements of a locally averaged quantity 
$q(x,t)$ that arises naturally in a standard finite volume discretisation of (\ref{eq:spde_constthetaeqn}).

To discretise the system, one first defines $q^i_t = q(x_i,t)$ around $N_d = 200$ grid points $x_i$ on a regular grid, separated by distances 
$\Delta x$, as
\begin{align}
  q_t^i =  \int_{x_i-\Delta x/2}^{x_i+\Delta x/2} u(x,t) \dif x .
\end{align}
The conservative (drift) term in (\ref{eq:spde_fulleqn}) reduces to
\begin{align}
 \int_{x_i-\Delta x/2}^{x_i+\Delta x/2} \nabla \cdot \left(\theta(x) \nabla u(x,t)\right) \dif x = 
 \theta_{i+1/2} \nabla u^{i+1/2}_t - \theta_{i-1/2} \nabla u^{i-1/2}_t ,
\end{align}
where $\theta_{i \pm 1/2} \equiv \theta(x_i + \Delta x/2)$, etc.  The following standard finite difference approximations
\begin{align}
  \nabla u^{i + 1/2}_t \simeq \frac{u^{i+1}_t - u^i_t}{\Delta x}, \quad u^i_t \simeq \Delta x^{-1} q^i_t
\end{align}
yield the $N_d$-dimensional SDE
\begin{align}
  \dif q_t^i  &= \theta \left(\frac{q^{i+1}_t - 2q^{i}_t + q^{i-1}_t}{\Delta x^2}\right) \dif t + \sigma^{1/2} \Delta x^{1/2} \dif W^i_t
\end{align}
for constant $\theta$, where $W_t^i$ are independent one-dimensional Brownian motions in time.  

Following recent results from \cite{sr:AR19} we consider the case of estimation of a constant $a = \theta$ value from measurements 
$\dif q^{\ast}_t$ at a fixed location/index $j^\ast \in \{1,\ldots,N_d\}$.  The data trajectory is thus given by 
\begin{align} \label{eq:spde_obs}
  \dif Y_t = \dif q^{\ast}_t + R^{1/2} \dif V_t
\end{align}
where $R^{1/2}$ is a scalar  and $V_t$ is a standard Brownian motion in one dimension.  We perform numerical experiments in which the initial state 
$q_0^i$ is set to zero for all indices $i$ and the prior on the unknown parameter $a = \theta$ is uniform over the interval $[0.2,1.8]$. 

The increment data is generated by first integrating (\ref{eq:spde_constthetaeqn}) forward in time from the known initial condition $q_i(0)=0$ for all $i$.  The equation is discretised in time using the Euler-Maruyama method.  It is known that $\Delta t < \theta \Delta x^2/2$ is required for stability of the 
Euler--Maruyama discretisation; we use the much smaller time step $\Delta t = \Delta x^2/80$.  The solution is sampled with this same time step, 
and increment measurements are approximated at time $t_n$ by setting the measurement noise level $R$ to zero in (\ref{eq:spde_obs}), resulting in
\begin{align}
\Delta Y_n = q^{\ast}_{n+1}-q^{\ast}_n.
\end{align}
Note that the associated model error in (\ref{eq:dynamics}) is given by $G = \sigma^{1/2} \Delta x^{1/2} I$ and the matrix $H$ in (\ref{eq:obs}) 
projects the vector of state increments onto a single component with index $j^\ast = N_d/2$. Simulations are performed over the time-interval $[0,20]$.
The results can be found in Figure \ref{fig:figure5}(a). We also compute the model evidence for a sequence of parameter values
$\theta \in \{0.2,0.3,\ldots,1.8\}$ based on a standard Kalman--Bucy filter \cite{sr:simon} 
for the associated linear state estimation problem. See Figure \ref{fig:figure5}(b).
Both approaches agree with the reference value $\theta = 1$.


\subsection{Discussion}

The presented results demonstrate that the proposed methodology can be applied to a broad range of continuous-time
state and parameter estimation problems with correlated measurement and model errors. Alternatively, one could have employed
standard SMC or MCMC methods utilising the modified observation model (\ref{eq:mod_forward}). However, such implementations require
the approximation of the additional $Q \nabla \log \pi_t$ term which is nontrivial if only samples from $\pi_t$ are available. Furthermore,
the limiting behaviour of such implementations in the limit $R\to 0$ and $H=I$ (pure parameter estimation problem) is unclear. The proposed
generalised ensemble Kalman--Bucy filter avoids these issues and is easy to implement. In fact, the only differences to the
standard ensemble Kalman--Bucy filter formulation of \cite{sr:br11} consist in the additional $QH^{\rm T}$ term in the Kalman 
gain and a correlation between the stochastic innovation process and the model error.

%
%
\section{Conclusions}
%
%

In this paper, we have derived McKean--Vlasov equations for combined state and parameter estimation from continuously observed
state increments. An approximate and robust implementation of these McKean-Vlassov equations in the form of a generalised ensemble
Kalman--Bucy filter has been provided and applied to a range of increasingly complex model systems. Future work will address
the treatment of temporally correlated measurement and model errors as well as a rigorous analysis of these McKean--Vlasov equations in a
multi-scale context and in the context of nonparametric drift estimation.


\vspace{6pt} 



\authorcontributions{Methodology, N.N. and S.R.; software, S.R. and P.R.; validation, N.N., S.R. and P.R.; writing--original draft preparation, 
N.N., S.R.; writing--review and editing, N.N., S.R. and P.R.}

\funding{This research has been partially funded by 
Deutsche Forschungsgemeinschaft (DFG) through grants 
CRC 1294 \lq Data Assimilation\rq \,(project A06) and CRC 1114 \lq Scaling Cascades\rq \,(project A02).}


\conflictsofinterest{The authors declare no conflict of interest. The funders had no role in the design of the study; in the collection, analyses, or 
interpretation of data; in the writing of the manuscript, or in the decision to publish the results.} 

%

\appendixtitles{no} 
\appendix
\section{The filtering equations for correlated noise}
In this appendix we outline a derivation of the Kushner-Stratonovich equation \eqref{eq:KS} for the signal-observation dynamics given by \eqref{eq:reformulation}. In fact, we only compute the evolution equation (termed  \emph{modified Zakai equation}) for the unnormalised filtering distribution $\rho_t[\phi] = \mathbb{E} \left[ l_t \phi(X_t) \vert Y_{[0,t]}\right]$, where the likelihood $l_t$ is given by
\begin{equation}
l_t \equiv l(Y_{[0,t]}|X_{[0,t]}) = \exp \left(\int_0^t f(X_s)^{\rm T} H^{\rm T} C^{-1} \,\mathrm{d}Y_s - \frac{1}{2} \int_0^t   f(X_s)^{\rm T} H^{\rm T} C^{-1}H f(X_s)\,{\rm d}s \right).
\end{equation}
Obtaining the Kushner-Stratonovich formulation is then standard, applying It\^o's formula to the Kallianpur-Striebel formula $\pi[\phi] = \rho_t[\phi] / \rho_t[\mathbf{1}]$, see \cite[Chapter 3]{sr:BC09}. The following result is in agreement with the corollaries 3.39 and 3.40 in \cite{sr:BC09}.

\begin{Lemma}
	The modified Zakai equation is given by
	\begin{equation}
	\rho_t[\phi] = \rho_0[\phi] + \int_0^t \rho_s[\mathcal{L} \phi] \, \mathrm{d}s + \int_0^t \rho_s\left[\phi f^{\rm T} 
	H^{\rm T} C^{-1}\right] \mathrm{d}Y_s + \int_0^t \rho_s \left[ \nabla \phi \right] QH^{\rm T} C^{-1}\, \mathrm{d}Y_s,
	\end{equation}
	where the generator $\mathcal{L}$ has been defined in \eqref{eq:generator}.
\end{Lemma}
\begin{proof}
	For convenience, let us define the process
	\begin{equation}
	\label{eq:obs process}
	M_t = \int_0^t f(X_s)^{\rm T} H^{\rm T} C^{-1} \mathrm{d}Y_s,
	\end{equation}
	where $Y_s$ satisfies (\ref{eq:obs2}).
	From $\langle Y \rangle_t = Ct$ we see that 
	\begin{equation}
	\langle M \rangle_t = \int_0^t   f(X_s)^{\rm T} H^{\rm T} C^{-1}H f(X_s)\,{\rm d}s,
	\end{equation}
	hence the likelihood takes the form 
	\begin{equation}
	l_t = \exp \left(M_t - \frac{1}{2} \langle M \rangle_t  \right)
	\end{equation}
	satisfying the SDE 
	\begin{equation}
	\label{eq:exp SDE2}
	\mathrm{d}l_t = l_t \, \mathrm{d}M_t.
	\end{equation}
	For an arbitrary smooth compactly supported test function $\phi$ It\^o's formula implies
	\begin{subequations}
		\label{eq:Ito for Zakai}
		\begin{align}
		l_t \phi(X_t) & = \phi(X_0) + \int_0^t \phi(X_s) \, \mathrm{d}l_s + \int_0^t l_s \nabla \phi(X_s) \cdot \mathrm{d}X_s \\
		& + \frac{1}{2} \sum_{i,j=1}^{N_x} \int_0^t l_s \partial_i \partial_j \phi(X_s) \, \mathrm{d}\langle X^i,X^j \rangle_s + \sum_{i=1}^{N_x}\int_0^t \partial_i \phi (X_s) \, \mathrm{d}\langle l, X^i\rangle_s,
		\end{align}
	\end{subequations}
	where $X_s$ satisfies (\ref{eq:dynamics2}).
	For the covariation process $\langle l, X\rangle_t$ we obtain
	\begin{equation}
	\langle l, X\rangle_t = l_t \langle M, X \rangle_t = l_t f(X_t)^{\rm T} H^{\rm T} C^{-1} HQ  t,
	\end{equation}
	using $\langle Y, X \rangle_t  = HQ  t$. Furthermore, $\langle X,X\rangle_t = Q t$, which follows from the
	definition of the stochastic contributions in (\ref{eq:dynamics2}).

	We now apply the conditional expectation to \eqref{eq:Ito for Zakai}. Noticing that
	\begin{equation}
	\int_0^t \phi(X_s) \, \mathrm{d}l_s = \int_0^t l_s \phi(X_s)f(X_s)^{\rm T} H^{\rm T} C^{-1} \, \mathrm{d}Y_s,
	\end{equation}
	the result follows from \eqref{eq:exp SDE2}.
\end{proof}



\reftitle{References}


\externalbibliography{yes}
\bibliography{bib_paper}



\end{document}